\documentclass[12pt]{amsart}
\usepackage{amssymb}

\usepackage{amsmath}

\theoremstyle{plain}

\numberwithin{equation}{section}
\input{tcilatex}

\begin{document}
\title[Free Modules, Semimodules and Free Lie Modules]{Automorphisms of
Categories of Free Modules, Free Semimodules, and Free Lie Modules}
\author{Yefim Katsov, Ruvim Lipyanski, and Boris Plotkin}
\address{\textit{Department of Mathematics and Computer Science}\\
\textit{Hanover College, Hanover, IN 47243--0890, USA}}
\email{\textit{katsov@hanover.edu}}
\address{\textit{Department of Mathematics}\\
\textit{Ben Gurion University, Beer Sheva, 84105, Israel}}
\email{\textit{lipyansk@cs.bgu.ac.il}}
\address{\textit{Institute of Mathematics}\\
\textit{The Hebrew University, Jerusalem, 91904, Israel}}
\email{\textit{borisov@macs.biu.ac.il}}
\date{April 30, 2005}
\dedicatory{Dedicated to the memory of Professor Saunders Mac Lane, a great
mathematician and man}

\begin{abstract}
In algebraic geometry over a variety of universal algebras $\Theta $, the
group $Aut(\Theta ^{0})$ of automorphisms of the category $\Theta ^{0}$ of
finitely generated free algebras of $\Theta $ is of great importance. In
this paper, semi-inner automorphisms are defined for the categories of free
(semi)modules and free Lie modules; then, under natural conditions on a
(semi)ring, it is shown that all automorphisms of those categories are
semi-inner. We thus prove that for a variety $_{R}\mathcal{M}$ of
semimodules over an IBN-semiring $R$ (an IBN-semiring is a semiring analog
of a ring with IBN), all automorphisms of $Aut(_{R}\mathcal{M}^{0})$ are
semi-inner. Therefore, for a wide range of rings, this solves Problem 12
left open in \cite{plotkin:slotuag}; in particular, for Artinian
(Noetherian, $PI$-) rings $R$, or a division semiring $R$, all automorphisms
of $Aut(_{R}\mathcal{M}^{0})$ are semi-inner.
\end{abstract}

\maketitle

\section{Introduction}

{\ {\ In the last decade, a{s an outcome of the major breakthrough in {%
Tarski's} striking conjecture about the elementarily equivalence {of free
groups} on more than one generator,} {a new algebraic area {{has appeared} }%
--- algebraic geometry over free group}s. Later {(see, for example, \cite%
{baummysrem:agog}), }algebraic geometry over free groups was extended to {%
algebraic geometry over arbitrary groups}, and even over groups with a fixed
group of constants. During the same time, in {{\cite{plotkin:voaaavcoav}, %
\cite{plotkin:snoagiua}, \cite{berzinsplts:agivoaw}, and \cite%
{plotkin:slotuag},}} {algebraic geometry} was, in fact, introduced{\ in }a
significantly more general setting --- namely, algebraic geometry over
algebras of an arbitrary variety of universal algebras (not only of group
varieties) --- that brought forth a fascinating {new} area of algebra known
as universal algebraic geometry. }}

As its name suggests, this area holds many surprising similarities to
classical algebraic geometry. However, in contrast to classical algebraic
geometry, which is closely connected with the ideal theory of finitely
generated polynomial algebras over fields, algebraic geometry in varieties
of universal algebras --- universal algebraic geometry --- is sharply
associated with congruence theories of finitely generated free algebras of
those varieties. Referring to \cite{plotkin:slotuag} and \cite%
{plotkin:snoagiua} for more details, this situation can be briefly described
as follows.

Let $\Theta $ be a variety of universal algebras, $X_{0}=\{x_{1},x_{2},%
\ldots ,x_{n},\ldots \}$ a fixed denumerable set of an infinite universe $%
\mathcal{U}$, and $\Theta ^{0}$ the full small ($\mathcal{U}$-small)
subcategory of\ finitely generated free algebras $F_{X},$ $X\ \subseteq 
\mathcal{U},|X|<$ $\infty $, of the variety $\Theta $. Then, for a fixed
algebra $G\in |\Theta |$ of the variety $\Theta $ and a free algebra $%
F_{X}\in |\Theta ^{0}|$, the set $Mor_{\Theta }(F_{X},G)$ of homomorphisms
from $F_{X}$ to $G$ is treated as an affine space over $G$ consisting of
points (homomorphisms) $\mu :F_{X}\longrightarrow G$. For any set of points $%
A\subseteq Mor_{\Theta }(F_{X},G)$ and any binary relation $T\subseteq
F_{X}\times F_{X}$ on $F_{X}$, the assignments $T\longmapsto T_{G}^{^{\prime
}}\overset{def}{=}\{\mu :F_{X}\longrightarrow G|$ $T\subseteq Ker\mu \}$ and 
$A\longmapsto A^{^{\prime }}\overset{def}{=}\cap _{\mu \in A}Ker\mu $ define
the Galois correspondence between binary relations (or systems of equations) 
$T\ $on $F_{X}$ and sets of points $A\ $of the space $Mor_{\Theta }(F_{X},G)$%
. A congruence $T$ on $F_{X}$ is said to be $G$-\textit{closed} if $%
T=A^{^{\prime }}$ for some point set $A\subseteq Mor_{\Theta }(F_{X},G)$; a
point set $A$ is called an \textit{algebraic variety} in the space $%
Mor_{\Theta }(F_{X},G)$ if $A=T_{G}^{^{\prime }}$ for some relation $T\ $on $%
F_{X}$. As usual, the Galois correspondence produces the closures: $%
A^{^{\prime \prime }}\overset{def}{=}(A^{^{\prime }})^{^{\prime }}$and $%
T^{^{\prime \prime }}\overset{def}{=}(T_{G}^{^{\prime }})^{^{\prime }}$.
Varying the free algebras $F_{X}\in |\Theta ^{0}|$, one comes up with the
category $K_{\Theta }(G)$ of algebraic varieties over $G$, which can be
regarded as an important geometric invariant of the algebra $G$ evaluating
``abilities'' of $G$ in solving systems of equations in free algebras $%
F_{X}\in |\Theta ^{0}|$. Thus, two algebras $G_{1}$, $G_{2}\in |\Theta |$
are said to be \textit{geometrically equivalent} iff $Cl_{G_{1}}(F_{X})%
\overset{def}{=}\{T$ $|$ $T$ is a $G_{1}$-closed relation on $F_{X}\}$ is
the same as $Cl_{G_{2}}(F_{X})\overset{def}{=}\{T$ $|$ $T$ is a $G_{2}$%
-closed relation on $F_{X}\}$ for any $F_{X}\in |\Theta ^{0}|$. The
definitions\ of $Cl_{G_{1}}(F_{X})$ and $Cl_{G_{2}}(F_{X})$ are naturally
extended to the (contravariant) functors $Cl_{G_{1}},Cl_{G_{2}}:\Theta
^{0}\longrightarrow Set$, and $G_{1}$ is geometrically equivalent\ to $G_{2}$
iff $Cl_{G_{1}}=Cl_{G_{2}}$.

One of the principal problems in algebraic geometry over a variety $\Theta $
involves studying interrelations between relations between algebras $G_{1}$
and $G_{2}$ and relations between $G_{1}$- and $G_{2}$-geometries over them;
in particular, relations between algebras $G_{1}$ and $G_{2}$ and the
isomorphic categories $K_{\Theta }(G_{1})$ and $K_{\Theta }(G_{2})$ (see,
for example, \cite{plotkin:voaaavcoav} and \cite{plotkin:slotuag}). Thus, it
is easy to show that the geometrical equivalence of algebras $G_{1}$ and $%
G_{2}$ implies that $K_{\Theta }(G_{1})$ and $K_{\Theta }(G_{2})$ are
isomorphic, \textit{i.e.}, $(Cl_{G_{1}}=Cl_{G_{2}}\Longrightarrow K_{\Theta
}(G_{1})$ $\cong $ $K_{\Theta }(G_{2}))$; but the converse implication is
not always true \cite[Theorem 6]{plotkin:voaaavcoav}. However, it is
important to know for which varieties $\Theta $ the equivalence $%
Cl_{G_{1}}=Cl_{G_{2}}\Longleftrightarrow K_{\Theta }(G_{1})$ $\cong $ $%
K_{\Theta }(G_{2})$ holds. (The reader may consult \cite{plotkin:voaaavcoav}%
, \cite{plotkin:snoagiua}, \cite{mashplts:aocfa}, and \cite{plotkin:slotuag}
for many interesting results regarding this problem.) The fundamental notion
of \textit{geometric similarity} of algebras, generalizing the notion of
geometric equivalence, proves to be crucial in all investigations concerning
this problem. In turn, the geometric similarity uses a more subtle
relationship between the functors $Cl_{G_{1}}$ and $Cl_{G_{2}}$ which, in
turn, is heavily based on automorphisms of the category $\Theta ^{0}$ (see,
for instance, \cite{plotkin:slotuag}). This fact has highly motivated a
sustained interest in studying automorphisms and autoequivalences of
categories $\Theta ^{0}$ for important varieties of universal algebras $%
\Theta $ (one may consult \cite{mashplts:aocfa}, \cite{plotkin:slotuag}, and %
\cite{masplotbe:acfla} for obtained results and open problems in this
direction).

This paper first of all concerns automorphisms of categories $\Theta ^{0}$
when the varieties $\Theta $ are quite important varieties of universal
algebras such as varieties of modules over rings, varieties of semimodules
over semirings, and varieties of Lie modules over Lie algebras. For all
these varieties, our main results are obtained as consequences of a general
original approach --- developed in the paper in the setting of semi-additive
categorical algebra and demonstrated in details for the most general case
--- for varieties of semimodules over semirings. Our special interest in
categories of semimodules over semirings is motivated, among other things,
by the following observations. Nowadays one may clearly notice a growing
interest in developing the algebraic theory of semirings and their numerous
connections with, and applications in, different branches of mathematics,
computer science, quantum physics, and many other areas of science (see, for
example, the recently published survey \cite{glazek:glosaa}). Also, as
algebraic objects, semirings certainly are the most natural generalization
of such (at first glance different) algebraic systems as rings and bounded
distributive lattices. Investigating semirings and their representations,
one should undoubtedly use methods and techniques of\ both ring and lattice
theory as well as diverse techniques and methods of categorical and
universal algebra. For these reasons, results originally obtained in
semiring and/or semimodule settings very often imply as immediate
consequences the corresponding results for well-established ``classical''
varieties such as rings, modules, semigroups, lattices, \textit{etc}.

For the reader's convenience, we include in Section 2 all subsequently
necessary notions and facts regarding automorphisms and autoequivalences of
categories. In this section, we also introduce some quite important concepts
--- an IBN-variety and a type of a variety --- for varieties of universal
algebras, which naturally extend the notions of an IBN-ring and a type of a
ring (see, for example, \cite{cohn:sroibp} and \cite{cohn:fratr}) to an
arbitrary variety of universal algebras, as well as establish some
interesting properties\ --- Propositions 2.9 and 2.10 --- connected with
these concepts and needed in the sequel.

In \cite[Theorem 2.4]{lipplt}, was proposed a computational proof that all
automorphisms of categories $_{R}\mathcal{M}^{0}$ of finitely generated free
modules over Noetherian rings are semi-inner. In fact, a similar proof can
be carried out for varieties of modules over any IBN-ring. Moreover, in the
setting of semi-additive categorical algebra in Section 3, we obtain the
main results of the paper --- Theorems 3.9, 3.10, and 3.13 --- essentially
improving and extending those results. Namely, we show that for (semi)rings
of almost all types, including IBN-(semi)rings, all automorphisms of
categories $_{R}\mathcal{M}^{0}$ of finitely generated free (semi)modules
are semi-inner; and therefore, positively resolve Problem 12 posted in \cite%
{plotkin:slotuag} for quite a wide spectrum of rings $R$. We also single out
Theorem 3.15 and Corollary 3.16, describing the group of outer automorphisms
of the (semi)module categories $_{R}\mathcal{M}^{0}$ in the most important
cases.

Then, continuing in the spirit of Sections 2 and 3, in Section 4 we obtain
Theorems 4.7, 4.10, 4.15, and Corollary 4.8, which are the Lie module
analogs of the corresponding main results of Section 3 for varieties of Lie
modules over Lie algebras.

We conclude the paper with Section 5, Appendix, providing a new
easy-to-prove version --- Proposition 5.1 and Theorem 5.3 --- of the quite
important reduction theorem (see \cite[Theorem 3.11]{mashplts:aocfa} and/or %
\cite[Theorem 3]{masplotbe:acfla}) that we use in the paper.

Finally, all notions and facts of categorical algebra, used here without any
comments, can be found in \cite{macl:cwm}; for notions and facts from
universal algebra, we refer to \cite{gratzr:ua}.

\section{Automorphisms of Categories}

Although it makes no relevant difference whether we base category theory on
the axiom of universes or on the G\"{o}del-Bernays theory of classes, in the
present context we have found it more reasonable to make use of universes
and require as an axiom that every set is an element of a universe. Thus,
from now on let $\mathcal{U}$ be a fixed universe containing the set $%
\mathbf{N}$ of natural numbers; and a category $\mathcal{C}$ is \textit{small%
}, more precisely $\mathcal{U}$-small, if the class $|\mathcal{C}|$ of all
objects of $\mathcal{C}$ is a $\mathcal{U}$-set.

As usual, an \textit{isomorphism} $\varphi :\mathcal{C}\longrightarrow 
\mathcal{D}$ of categories is a functor $\varphi $ from $\mathcal{C}$ to $%
\mathcal{D}$ which is a bijection, both on objects and morphisms; or
alternatively, it can be defined as a functor $\varphi $ from $\mathcal{C}$
to $\mathcal{D}$ for which there exists a functor $\psi $ from $\mathcal{D}$
to $\mathcal{C}$ such that $\psi \varphi =1_{\mathcal{C}}$ and $\varphi \psi
=1_{\mathcal{D}}$. We will also need a weaker notion than that of
isomorphism --- the notion of equivalence of two categories. Namely, an 
\textit{equivalence} between categories $\mathcal{C}$ and $\mathcal{D}$ is
defined to be a pair of functors $\varphi :\mathcal{C}\longrightarrow 
\mathcal{D}$, $\psi :\mathcal{D}\longrightarrow \mathcal{C}$ together with
natural isomorphisms $\psi \varphi \cong 1_{\mathcal{C}}$ and $\varphi \psi
\cong 1_{\mathcal{D}}$. It is easy to see (also see, for example, %
\cite[Propositions 16.3.2, 16.3.4, and 16.3.6]{schu:cat}) that in any
category $\mathcal{C}$ there exists a full subcategory\ $\mathcal{C}_{Sk}$
of $\mathcal{C}$, called a \textit{skeleton} of $\mathcal{C}$, such that
each object of $\mathcal{C}$ is isomorphic (in $\mathcal{C}$) to exactly one
object in $\mathcal{C}_{Sk}$; and therefore, the inclusion $I:\mathcal{C}%
_{Sk}\longrightarrow \mathcal{C}$ defines an equivalence of categories.

Let $\mathcal{C}$ be a small category. Then, all endofunctors $\varphi :%
\mathcal{C}\longrightarrow \mathcal{C}$ with respect to composition form the
semigroup $End$ $(\mathcal{C})$ of endomorphisms of the category $\mathcal{C}
$, which contains the subgroup $Aut$ $(\mathcal{C})$ of automorphisms of the
category $\mathcal{C}$. We will distinguish the following classes of
automorphisms of $\mathcal{C}$.\medskip\ \ 

\noindent \textbf{Definition 2.1. }An automorphism $\varphi :\mathcal{C}%
\longrightarrow \mathcal{C}$ is \textit{equinumerous} if $\varphi (A)\cong A$
for any object $A\in |\mathcal{C}|$; $\varphi $ is \textit{stable} if $%
\varphi (A)=A$ for any object $A\in |\mathcal{C}|$; and $\varphi $ is 
\textit{inner} if $\varphi $ and $1_{\mathcal{C}}$ are naturally
isomorphic.\medskip

The following observation is obvious.\medskip

\noindent \textbf{Proposition 2.2. }\textit{Within the group }$Aut$ $(%
\mathcal{C})$\textit{,} \textit{the collections} $EqnAut$ $(\mathcal{C})$, $%
StAut$ $(\mathcal{C})$\textit{, and} $Int$ $(\mathcal{C})$ \textit{of
equinumerous, stable, and inner automorphisms, respectively, form normal
subgroups. Also,} $Int$ $(\mathcal{C})\subseteq EqnAut$ $(\mathcal{C})$ 
\textit{and} $StAut$ $(\mathcal{C})\subseteq EqnAut$ $(\mathcal{C})$\textit{%
.\medskip }

\noindent \textbf{Definition 2.3. }A group $Out$ $(\mathcal{C})$ of outer
automorphisms of the category $\mathcal{C}$ is defined as the quotient group 
$Aut$ $(\mathcal{C})/Int$ $(\mathcal{C})$, and $\mathcal{C}$ is called 
\textit{perfect} if this group is trivial.\medskip

Let $\mathcal{C}_{Sk}$ be a skeleton of $\mathcal{C}$. Then, for each object
object $A\in |\mathcal{C}|$, there exists a unique object $\overline{A}\in |%
\mathcal{C}_{Sk}|$ isomorphic to $A$; and let $\{i_{A}:A\longrightarrow 
\overline{A\text{ }}|$ $A\in |\mathcal{C}|$, and $i_{A}=1_{A}$ if $A=%
\overline{A}\}$ be a fixed set of isomorphisms of the category $\mathcal{C}$%
. The following observations will prove to be useful.\medskip

\noindent \textbf{Lemma 2.4.} \textit{The assignment} $\varphi \longmapsto 
\overline{\varphi }\overset{def}{=}\varphi |_{\mathcal{C}_{Sk}}:\mathcal{C}%
_{Sk}\longrightarrow \mathcal{C}_{Sk}$\textit{\ --- the restriction of a
stable automorphism }$\varphi \in StAut$\textit{\ }$(\mathcal{C})$\textit{\
to the subcategory }$\mathcal{C}_{Sk}$\textit{\ --- defines a group
homomorphism }

\noindent $|_{\mathcal{C}_{Sk}}:$\textit{\ }$StAut$\textit{\ }$(\mathcal{C}%
)\longrightarrow StAut$\textit{\ }$(\mathcal{C}_{Sk})$\textit{.\medskip }

\noindent \textbf{Proof}. Indeed, as obviously $\overline{\varphi }$ is, in
fact, the composite $I\varphi I=\varphi I:\mathcal{C}_{Sk}\longrightarrow 
\mathcal{C}_{Sk}$ for any $\varphi \in StAut$ $(\mathcal{C})$, one has $%
\overline{\psi }\overline{\varphi }=$ $I\psi II\varphi I=I\psi I\varphi
I=I\psi \varphi I=\overline{\psi \varphi }$ for any $\psi ,\varphi \in StAut$
$(\mathcal{C})$.\textit{\ \ \ \ \ \ }$_{\square }$\smallskip \medskip

\noindent \textbf{Lemma 2.5. }\textit{Let} $\overline{\varphi }\in StAut$%
\textit{\ }$(\mathcal{C}_{Sk})$\textit{. Then the assignment} $\overline{%
\varphi }\longmapsto \overline{\varphi }^{i}:\mathcal{C}\longrightarrow 
\mathcal{C}$\textit{, where} $\overline{\varphi }^{i}(f)\overset{def}{=}%
i_{B}^{-1}\overline{\varphi }(i_{B}f$ $i_{A}^{-1})$ $i_{A}$ \textit{for any
morphism} $f\in Mor_{\mathcal{C}}$\textit{\ }$(A,B)$\textit{, defines a
group monomorphism} $^{i}:StAut$\textit{\ }$(\mathcal{C}_{Sk})%
\longrightarrow StAut$\textit{\ }$(\mathcal{C})$\textit{.}\medskip

\noindent \textbf{Proof}. First, if $g\in Mor_{\mathcal{C}}$\textit{\ }$%
(B,C) $, then $\overline{\varphi }^{i}(gf)=i_{C}^{-1}\overline{\varphi }%
(i_{C}$ $gf $ $i_{A}^{-1})$ $i_{A}=i_{C}^{-1}\overline{\varphi }(i_{C}$ $g$ $%
i_{B}^{-1}i_{B}f$ $i_{A}^{-1})$ $i_{A}=i_{C}^{-1}\overline{\varphi }(i_{C}$ $%
g$ $i_{B}^{-1})\overline{\varphi }(i_{B}f$ $i_{A}^{-1})$ $i_{A}=i_{C}^{-1}%
\overline{\varphi }(i_{C}$ $g$ $i_{B}^{-1})$ $i_{B}i_{B}^{-1}\overline{%
\varphi }(i_{B}f$ $i_{A}^{-1})$ $i_{A}=\overline{\varphi }^{i}(g)\overline{%
\varphi }^{i}(f)$. Hence, indeed $\overline{\varphi }^{i}\in StAut$\textit{\ 
}$(\mathcal{C})$.

Now, if $\overline{\varphi }$, $\overline{\psi }\in StAut$\textit{\ }$(%
\mathcal{C}_{Sk})$, then $(\overline{\varphi }\overline{\psi }%
)^{i}(f)=i_{B}^{-1}\overline{\varphi }\overline{\psi }(i_{B}f$ $i_{A}^{-1})$ 
$i_{A}=i_{B}^{-1}\overline{\psi }(\overline{\varphi }(i_{B}f$ $i_{A}^{-1}))$ 
$i_{A}=i_{B}^{-1}\overline{\psi }(i_{B}i_{B}^{-1}\overline{\varphi }(i_{B}f$ 
$i_{A}^{-1})i_{A}i_{A}^{-1})$ $i_{A}=$

\noindent $i_{B}^{-1}\overline{\psi }(i_{B}\overline{\varphi }%
^{i}(f)i_{A}^{-1})$ $i_{A}=\overline{\psi }^{i}(\overline{\varphi }^{i}(f))=%
\overline{\varphi }^{i}\overline{\psi }^{i}(f)$. Therefore, $^{i}:StAut$%
\textit{\ }$(\mathcal{C}_{Sk})\longrightarrow StAut$\textit{\ }$(\mathcal{C}%
) $ is a group homomorphism, which is obviously mono.\textit{\ \ \ \ \ \ }$%
_{\square }\medskip $

\noindent \textbf{Proposition 2.6. }\textit{For any} $\varphi \in StAut$%
\textit{\ }$(\mathcal{C})$\textit{, the automorphisms }$\varphi $ \textit{and%
} $(|_{\mathcal{C}_{Sk}}(\varphi ))^{i}=\overline{\varphi }^{i}\in StAut$%
\textit{\ }$(\mathcal{C})$ \textit{are naturally isomorphic, i.e.,} $%
\overline{\varphi }^{i}\cong $ $\varphi $ \textit{in the category} $\mathcal{%
F}(\mathcal{C},\mathcal{C})$ \textit{of endofunctors on} $\mathcal{C}$%
\textit{.}\medskip

\noindent \textbf{Proof}. Indeed, for any $f\in Mor_{\mathcal{C}}$\textit{\ }%
$(A,B)$ we have $\overline{\varphi }^{i}(f)=i_{B}^{-1}\overline{\varphi }%
(i_{B}f$ $i_{A}^{-1})i_{A}=i_{B}^{-1}\varphi (i_{B}f$ $%
i_{A}^{-1})i_{A}=i_{B}^{-1}\varphi (i_{B})\varphi (f)\varphi
(i_{A}^{-1})i_{A}$, and, hence,

\noindent $\varphi (f)\varphi (i_{A}^{-1})i_{A}=\varphi (i_{B}^{-1})i_{B}%
\overline{\varphi }^{i}(f)$. Thus, the set $\{\varphi
(i_{A}^{-1})i_{A}:A\longrightarrow A$ $|$ $A\in |\mathcal{C}|\}$ of
isomorphisms sets a natural isomorphism $\overline{\varphi }^{i}\overset{%
\bullet }{\longrightarrow }$ $\varphi $.\textit{\ \ \ \ \ \ }$_{\square
}\medskip $

\noindent \textbf{Proposition 2.7. }\textit{For any equinumerous automorphism%
} $\varphi \in EqnAut$ $(\mathcal{C})$ \textit{there exist a stable
automorphism} $\varphi _{S}\in StAut$\textit{\ }$(\mathcal{C})$ \textit{and
an inner automorphism} $\varphi _{I}\in Int$ $(\mathcal{C})$ \textit{such
that} $\varphi =\varphi _{S}\varphi _{I}$\textit{.}\medskip

\noindent \textbf{Proof}. Let $\varphi \in EqnAut$ $(\mathcal{C})$, and $%
\{i_{A}:A\longrightarrow \varphi (A)$ $|$ $A\in |\mathcal{C}|\}$ be a fixed
set of isomorphisms of the category $\mathcal{C}$. Then, define $\varphi
_{S}(A)\overset{def}{=}A$ for any $A\in |\mathcal{C}|$, and $\varphi _{S}(f)%
\overset{def}{=}$ $i_{B}^{-1}\varphi (f)i_{A}$ for any morphism $f\in Mor_{%
\mathcal{C}}$\textit{\ }$(A,B)$; and $\varphi _{I}(A)\overset{def}{=}\varphi
(A)$ for any $A\in |\mathcal{C}|$, and $\varphi _{I}(f)\overset{def}{=}%
i_{B}f $ $i_{A}^{-1}$ for any morphism $f\in Mor_{\mathcal{C}}$\textit{\ }$%
(A,B)$. From these definitions one can immediately see that $\varphi _{S}\in
StAut$\textit{\ }$(\mathcal{C})$, $\varphi _{I}\in Int$ $(\mathcal{C})$, and 
$\varphi _{S}\varphi _{I}(f)=\varphi _{I}(\varphi _{S}(f))=\varphi
_{I}(i_{B}^{-1}\varphi (f)i_{A})=\varphi (f)$, \textit{i.e}., $\varphi
=\varphi _{S}\varphi _{I}$\textit{.\ \ \ \ \ \ }$_{\square }\medskip $

In this paper, our main interest concerns the following situation. Let $%
\Theta $ be a variety of universal algebras, $X_{0}=\{x_{1},x_{2},\ldots
,x_{n},\ldots \}$ $\subseteq $ $\mathcal{U}$ a fixed denumerable set of the
infinite universe $\mathcal{U}$, and $\Theta ^{0}$ the full subcategory of\
finitely generated free algebras $F_{X},$ $X\ \subseteq \mathcal{U},|X|<$ $%
\infty $, of the variety $\Theta $. Then, let $F_{n},n\in $ $\mathbf{N}$,
stand for a free algebra generated by the elements $x_{1},x_{2},\ldots
,x_{n} $ of $X_{0}$, and consider the full subcategory $\Theta _{Sk}^{0}$ of
the category $\Theta ^{0}$ defined by the algebras $F_{n},n\in $ $\mathbf{N}$%
. Note: if the category $\Theta $ contains a fixed one-element zero object $%
F_{0}\in |\Theta |$, we will consider $F_{0}$ to be a free algebra with the
empty set of generators and, therefore, include it in $|\Theta _{Sk}^{0}|$.
\medskip

\noindent \textbf{Definition 2.8.} A variety $\Theta $ is said to have IBN
(``\textit{Invariant Basis Number}''), or to be an IBN-variety, if for any
natural numbers $n$, $m$ an isomorphism between free algebras with $n$ and $%
m $ free generators implies that $n=m$. \medskip

Thus, if the variety $\Theta $ considered above has IBN, then the
subcategory $\Theta _{Sk}^{0}$ is obviously a skeleton of the category $%
\Theta ^{0}$, and we have the following fact.\medskip

\noindent \textbf{Proposition 2.9.} \textit{Let} $\Theta $ \textit{be an
IBN-variety. Then all autoequivalences }$\varphi :\Theta ^{0}\longrightarrow
\Theta ^{0}$ \textit{and} $\overline{\varphi }:\Theta
_{Sk}^{0}\longrightarrow \Theta _{Sk}^{0}$ \textit{of the categories} $%
\Theta ^{0}$ \textit{and} $\Theta _{Sk}^{0}$\textit{, respectively, are
equinumerous, i.e., }$\varphi (A)\cong A$ \textit{and} $\overline{\varphi }(%
\overline{A})=\overline{A}$ \textit{for any objects} $A\in |\Theta ^{0}|$ 
\textit{and} $\overline{A}\in |\Theta _{Sk}^{0}|$\textit{. In particular,} $%
Aut$ $(\Theta ^{0})=EqnAut$ $(\Theta ^{0})$ \textit{and} $Aut$ $(\Theta
_{Sk}^{0})=EqnAut$ $(\Theta _{Sk}^{0})=StAut$\textit{\ }$(\Theta _{Sk}^{0})$%
\textit{.\medskip }

\noindent \textbf{Proof}. Let $\varphi :\Theta ^{0}\longrightarrow \Theta
^{0}$ be an equivalence of $\Theta ^{0}$. Then, there exists an endofunctor $%
\psi :\Theta ^{0}\longrightarrow \Theta ^{0}$ together with natural
isomorphisms $\psi \varphi \cong 1_{\Theta ^{0}}$ and $\varphi \psi \cong
1_{\Theta ^{0}}$. By \cite[Theorem IV.4.1]{macl:cwm} (see also \cite[Remark
16.5.9]{schu:cat}), the functor $\psi $ is the left and right adjoint of the
functor $\varphi $,\textit{\ i.e.}, $\psi \dashv \varphi $ and $\varphi
\dashv \psi $. Therefore, by \cite[Theorem V.5.1]{macl:cwm} (see also %
\cite[Proposition 16.2.4]{schu:cat}), $\varphi $ and $\psi $ preserve all
limits and colimits, and (in particular) all products and coproducts. From
the latter, taking into account that any algebra $F_{X}\in |\Theta ^{0}|$ is
a coproduct of $X$ copies of the free algebra $F_{1}\in |\Theta ^{0}|$ by %
\cite[Corollary 12.11]{mal:as} (also see \cite{sikorski:poaa}), one
immediately obtains the statement for the equivalence $\varphi :\Theta
^{0}\longrightarrow \Theta ^{0}$.

The same arguments work well for an equivalence $\overline{\varphi }:\Theta
_{Sk}^{0}\longrightarrow \Theta _{Sk}^{0}$, too.\textit{\ \ \ \ \ \ }$%
_{\square }\medskip $

Thus, without loss of generality, by \cite[Corollary 12.11]{mal:as} we may
consider the free algebras $F_{n}\in |\Theta ^{0}|$ to be coproducts of $n$
copies of $F_{1}\in |\Theta ^{0}|$, \textit{i.e.}, $F_{n}=\amalg
_{i=1}^{n}F_{1i}$, $F_{1i}=F_{1}$, $i=1,\ldots ,n$, $n\in $ $\mathbf{N}$.
Then, if $\Theta $ is not an IBN-variety, we have $F_{n}\cong F_{m}$ for
some $n,m\in $ $\mathbf{N}$, $n\neq m$, which is apparently the same as to
say that there exist two natural numbers $n,h\in $ $\mathbf{N}$ such that $%
F_{n}\cong F_{n+h}$. Following \cite{cohn:fratr} (see also \cite{cohn:sroibp}%
), the first pair $(n,h)$ in the usual lexicographic ordering satisfying $%
F_{n}\cong F_{n+h}$ is called the \textit{type} of $F_{1}$, or of $\Theta $.
From this, it is obvious that up to isomorphism there are only $n+h-1$
different finitely generated free algebras in $\Theta $ of the type $(n,h)$,
namely $F_{1},F_{2},\ldots ,F_{n}\ldots ,F_{n+h-1}$. Thus, the full
subcategory $\Theta _{Sk}^{0}$ of the category $\Theta ^{0}$, consisting of $%
F_{1},F_{2},\ldots ,F_{n}\ldots ,F_{n+h-1}$ and perhaps plus a zero object $%
F_{0}$, is a skeleton of $\Theta ^{0}$. One can easily note (see also, %
\cite[Section 12.2]{mal:as}) that with respect to the coproduct operation
the free algebras $F_{1},F_{2},\ldots ,F_{n}\ldots ,F_{n+h-1}$ form a
monogenic semigroup $M(n,h)$ with index $n$ and period $h$ and generated by $%
F_{1}$. In this context, we have the following useful observation:\medskip

\noindent \textbf{Proposition 2.10. } \textit{If }$\Theta $\textit{\ has
type }$(1,1)$\textit{, or} $(1,2)$, \textit{or} $(n,h)$ \textit{with} $n>1$,%
\textit{\ then any automorphism }$\varphi :\Theta ^{0}\longrightarrow \Theta
^{0}$\textit{\ is equinumerous.} \medskip

\noindent \textbf{Proof}. As $\Theta _{Sk}^{0}$ is a skeleton of $\Theta
^{0} $, an autoequivalence $\varphi :\Theta ^{0}\longrightarrow \Theta ^{0}$
induces the autoequivalence $\overline{\varphi }:\Theta
_{Sk}^{0}\longrightarrow \Theta _{Sk}^{0}$ defined as $\overline{\varphi }%
\overset{def}{=}I^{-1}\varphi I$, where $I^{-1}$ is an equivalence-inverse
to the inclusion $I$ (see, for example, \cite[Definition 16.2.1]{schu:cat}).
As in Proposition 2.9, $\overline{\varphi }$ preserves coproducts, and
hence, the generator of the monogenic semigroup $M(n,h)$. Therefore, $%
\varphi (F_{1})\cong F_{1}$ for an automorphism $\varphi $, and $\varphi $
is equinumerous.\textit{\ \ \ \ \ \ }$_{\square }\medskip $

We end this section with the following quite natural and interesting open
problem.\medskip

\noindent \textbf{Problem 1. }Does there exist a variety $\Theta $ of the
type $(1,h)$ with $h\geq 3$ such that not all automorphisms $\varphi :\Theta
^{0}\longrightarrow \Theta ^{0}$\textit{\ }are equinumerous? (With respect
to this question, we have a conjecture: There exists a ring $R$ such that
the variety $_{R}Mod$\textit{\ }of left modules over $R$ has a type $(1,h)$
with $h\geq 3$, and not all automorphisms $\varphi :$ $_{R}Mod^{0}%
\longrightarrow $ $_{R}Mod^{0}$\textit{\ }are equinumerous.)\textit{\ }

\section{Automorphisms of Categories of Semimodules}

Recall \cite{golan:sata} that a \emph{semiring\/} is an algebra $(R,+,\cdot
,0,1)$ such that the following conditions are satisfied:

(1) $(R,+,0)$ is a commutative monoid with identity element $0$;

(2) $(R,\cdot ,1)$ is a monoid with identity element $1$;

(3) Multiplication distributes over addition from either side;

(4) $0r=0=r0$ for all $r\in R$.

As usual, a \emph{left\/} $R$-\emph{semimodule} over the semiring $R$ is a
commutative monoid $(M,+,0_{M})$ together with a scalar multiplication $%
(r,m)\mapsto rm$ from $R\times M$ to $M$ which satisfies the following
identities for all $r,r^{^{\prime }}\in R$ and $m,m^{^{\prime }}\in M$:

(1) $(rr^{^{\prime}})m=r(r^{^{\prime}}m)$;

(2) $r(m+m^{^{\prime}})=rm+rm^{^{\prime}}$;

(3) $(r+r^{^{\prime }})m=rm+rm^{^{\prime }}$;

(4) $1m=m$;

(5) $r0_{M}=0_{M}=0m$.

\emph{Right semimodules\/} over $R$ and homomorphisms between semimodules
are defined in the standard manner. And, from now on, let $\mathcal{M}$ be
the variety of commutative monoids, and $\mathcal{M}_{R}$ and $_{R}\mathcal{M%
}$ denote the categories of right and left semimodules, respectively, over a
semiring $R$. As usual (see, for example, \cite[Chapter 17]{golan:sata}), if 
$R$ is a semiring, then in the category $_{R}\mathcal{M}$, a \textit{free}
(left) semimodule $\sum_{i\in I}R_{i},R_{i}\cong $ $_{R}R$, $i\in I$, with a
basis set $I$ is a direct sum (a coproduct) of $I$-th copies of $_{R}R$; and
free right semimodules are defined similarly.

Following \cite[Section 1.5]{schu:cat}, we say that a category $\mathcal{C}$
is \textit{semi-additive} if $Mor_{\mathcal{C}}$\textit{\ }$(A,B)$ is a
commutative monoid for any two objects $A,B\in |\mathcal{C}|$, and the
composition of morphisms is distributive on both sides and is compatible
with $0$-elements. It is clear that $\mathcal{M}_{R}$ and $_{R}\mathcal{M}$
are semi-additive categories, as well as are any full subcategories of these
categories; in particular, any full subcategories of free right and left
semimodules of $\mathcal{M}_{R}$ and $_{R}\mathcal{M}$ are always
semi-additive. By \cite[Proposition 12.2.5 and Convention 12.2.6]{schu:cat},
a semi-additive category with a zero object has finite products iff it has
finite coproducts. If this is the case, then finite coproducts are also
products, and we then talk about \textit{biproducts} and use the notation $%
\oplus A_{i}\overset{\pi _{i}}{\underset{\mu _{i}}{\rightleftarrows }}A_{i}$%
, where injections\ $\mu _{i}$ and projections $\pi _{i}$ are subjected to
the conditions $\sum_{i}\pi _{i}\mu _{i}=1_{\oplus A_{i}}$, $\mu _{i}\pi
_{i}=1_{A_{i}}$, and $\mu _{i}\pi _{j}=0$ if $i\neq j$.

A semiring $R$ is a left (right) \textit{IBN-semiring} if $_{R}\mathcal{M}$ (%
$\mathcal{M}_{R}$) is an IBN-variety. The following observation shows that
we can speak of an ``IBN-semiring'' without specifying ``right'' or
``left.'' \medskip

\noindent \textbf{Proposition 3.1.} \textit{A semiring }$R$ \textit{is a
left IBN-semiring iff it is a right IBN-semiring.\medskip }

\noindent \textbf{Proof}. $\Rightarrow $. Suppose $R$ is a left
IBN-semiring, and $\oplus _{i=1}^{n}R_{i}$ and $\oplus _{j=1}^{m}R_{j}$,
where $R_{i}\cong $ $R_{R}\cong R_{j}$ for all $i$, $j$, are isomorphic for
some $n,m\in \mathbf{N}$. Then, from \cite[Proposition 3.8]{katsov:tpies}
and \cite[Theorem 3.3]{katsov:thhs}, one readily has $\oplus _{i=1}^{n}$ $%
_{R}R_{i}\cong (\oplus _{i=1}^{n}$ $_{R}R_{i})\otimes _{R}R\cong (\oplus
_{j=1}^{m}$ $_{R}R_{j})\otimes _{R}R\cong \oplus _{j=1}^{m}$ $_{R}R_{j}$,
and, hence, $n=m$, and $R$ is a right IBN-semiring.

$\Leftarrow $. This is proved in a similar fashion.\textit{\ \ \ \ \ \ }$%
_{\square }\medskip $

From now on, we assume that $R$ is an IBN-semiring; $_{R}\mathcal{M}^{0}$ is
the full subcategory of\ finitely generated free (left) $R$-semimodules $%
F_{X},$ $X\ \subseteq \mathcal{U},|X|<$ $\infty $, of the variety $_{R}%
\mathcal{M}$; and $_{R}\mathcal{M}_{Sk}^{0}$denotes the full subcategory of
free semimodules $F_{n}$, generated by the elements $x_{1},x_{2},\ldots
,x_{n}$ of $X_{0}=\{x_{1},x_{2},\ldots ,x_{n},\ldots \}$ $\subseteq $ $%
\mathcal{U},$ $n\in $ $\mathbf{N}$. Apparently both $_{R}\mathcal{M}^{0}$
and $_{R}\mathcal{M}_{Sk}^{0}$ are semi-additive categories with biproducts.
Also it is clear that we may, without loss of generality, accept that $%
F_{n}= $ $\oplus _{i=1}^{n}R_{i}\overset{\pi _{i}}{\underset{\mu _{i}}{%
\rightleftarrows }}R_{i}$, where $R_{i}=$ $_{R}R$ with the canonical
injections\ $\mu _{i}$ and projections $\pi _{i}$ for any $i=1,\ldots ,n$
and $n\in $ $\mathbf{N}$. As $R$ is an IBN-semiring, $_{R}\mathcal{M}%
_{Sk}^{0}$ is a skeleton of the category $_{R}\mathcal{M}$, and, by
Proposition 2.9, all automorphisms of $_{R}\mathcal{M}_{Sk}^{0}$ are stable.

Now, let $\sigma :R\longrightarrow R$ be a semiring automorphism. For any $%
n\in $ $\mathbf{N}$, there exists the monoid automorphism $\sigma _{n}:$ $%
\oplus _{i=1}^{n}R_{i}\longrightarrow \oplus _{i=1}^{n}R_{i}$ defined by the
assignment $\oplus _{i=1}^{n}a_{i}\longmapsto $ $\oplus
_{i=1}^{n}(a_{i})^{\sigma }$. Then, one can readily verify that there is the
automorphism $\varphi ^{\sigma }:$ $_{R}\mathcal{M}_{Sk}^{0}\longrightarrow $
$_{R}\mathcal{M}_{Sk}^{0}$ defined by the assignments $\varphi ^{\sigma
}(f)(\oplus _{i=1}^{n}a_{i})\overset{def}{=}$ $\sigma _{m}(f((\oplus
_{i=1}^{n}a_{i}^{\sigma ^{-1}}))$ for any $\oplus _{i=1}^{n}a_{i}\in \oplus
_{i=1}^{n}R_{i}$ and $f\in Mor_{_{R}\mathcal{M}_{Sk}^{0}}(\oplus
_{i=1}^{n}R_{i},\oplus _{j=1}^{m}R_{j})$. Therefore, the following diagram
is commutative: 
\begin{equation*}
\begin{array}{ccc}
F_{n}=\oplus _{i=1}^{n}R_{i} & \overset{\sigma _{n}}{\longrightarrow } & 
\oplus _{i=1}^{n}R_{i}=F_{n} \\ 
f\downarrow &  & \downarrow \varphi ^{\sigma }(f) \\ 
F_{m}=\oplus _{j=1}^{m}R_{j} & \underset{\sigma _{m}}{\longrightarrow } & 
\oplus _{j=1}^{m}R_{j}=F_{m}%
\end{array}
\text{ .\ \ \ \ \ \ \ \ \ \ \ \ \ \ \ \ \ \ \ \ \ \ \ \ \ (1)}
\end{equation*}

\noindent \textbf{Definition 3.2.} An automorphism $\varphi \in Aut$ $(_{R}%
\mathcal{M}_{Sk}^{0})=StAut$\textit{\ }$(_{R}\mathcal{M}_{Sk}^{0})$ is
called \textit{skew-inner} if $\varphi =\varphi ^{\sigma }$ for some
semiring automorphism $\sigma :R\longrightarrow R$.\medskip

\noindent \textbf{Proposition 3.3.} \ \textit{Every} $\varphi \in Aut$ $(_{R}%
\mathcal{M}_{Sk}^{0})$\textit{\ that is constant on the canonical injections 
}$\mu _{i}:R_{i}\longrightarrow \oplus _{i=1}^{n}R_{i}$\textit{, i.e.,} $%
\varphi (\mu _{i})=\mu _{i}$ \textit{for all} $i=1,\ldots ,n$ \textit{and} $%
n\in $ $\mathbf{N}$\textit{,} \textit{is skew-inner.}\textbf{\ \medskip }

\noindent \textbf{Proof}. First, note that if $f:F_{n}=\oplus
_{i=1}^{n}R_{i}\longrightarrow \oplus _{j=1}^{m}R_{j}=F_{m}$ is a
homomorphism of left $R$-semimodules, and $\oplus _{i=1}^{n}R_{i}\overset{%
\pi _{i}}{\underset{\mu _{i}}{\rightleftarrows }}R_{i}$ and $\oplus
_{j=1}^{m}R_{j}\overset{\rho _{j}}{\underset{\lambda _{j}}{\rightleftarrows }%
}R_{j}$ are the canonical injections and projections corresponding to the
biproducts, then it is clear (also see, for example, \cite[Section 12.2.]%
{schu:cat}) that $f=\oplus _{i}((\times )_{j}(\mu _{i}f\rho _{j}))$. Hence, $%
\varphi (f)=\oplus _{i}((\times )_{j}\varphi (\mu _{i}f\rho _{j}))$ for any $%
\varphi \in Aut$ $(_{R}\mathcal{M}_{Sk}^{0})$ since, as was mentioned in the
proof of Proposition 2.9, $\varphi $ preserves all colimits and limits, and,
in particular, biproducts. Thus, $\varphi (f)$ is completely defined by $%
\varphi (\mu _{i}f\rho _{j}):$ $_{R}R=R_{i}\longrightarrow R_{j}=$ $_{R}R$
for $\mu _{i}f\rho _{j}:$ $_{R}R=R_{i}\longrightarrow R_{j}=$ $_{R}R$, $%
i=1,\ldots ,n$, $j=1,\ldots ,m$. In turn, $\varphi (\mu _{i}f\rho _{j})$ is
defined by $\varphi _{R}$, the action of the functor $\varphi $ on the
semiring $End(_{R}R)=Mor_{_{R}\mathcal{M}_{Sk}^{0}}(_{R}R,_{R}R)$ of
endomorphisms of the regular semimodule $_{R}R$. Then, agreeing to write
endomorphisms of $_{R}R$ on the right of the elements they act on, one can
easily see that actions of endomorphisms actually coincide with
multiplications of elements of $_{R}R$ on the right by elements of the
semiring $R$, and, therefore, $End(_{R}R)=$ $R$; and $\varphi _{R}$ is a
monoid automorphism of the multiplicative reduct of $R$. Furthermore, since
the category $_{R}\mathcal{M}_{Sk}^{0}$ has a zero object and finite
products, by \cite[Proposition 12.2.7]{schu:cat} the functor $\varphi $ is
additive; and, hence, $\varphi _{R}$ is also a monoid automorphism of the
additive reduct of $R$, \textit{i.e.}, $\varphi _{R}$ is just a semiring
automorphism of $R$ that we denote by $\sigma $.

Secondly, because $\varphi $ is an additive\ functor preserving biproducts,
from $\sum_{i}\pi _{i}\mu _{i}=1_{\oplus _{i=1}^{n}R_{i}}$, $\mu _{i}\pi
_{i}=1_{R_{i}}$, and $\mu _{i}\pi _{j}=0$ if $i\neq j$, one has $\varphi
(\sum_{i}\pi _{i}\mu _{i})=\sum_{i}\varphi (\pi _{i}\mu
_{i})=\sum_{i}\varphi (\pi _{i})\varphi (\mu _{i})=\sum_{i}\varphi (\pi
_{i})\mu _{i}=1_{\oplus _{i=1}^{n}R_{i}}$, and, therefore, $\varphi (\pi
_{i})=(\sum_{i}\varphi (\pi _{i})\mu _{i})\pi _{i}=\pi _{i}$ for all $%
i=1,\ldots ,n$ and $n\in $ $\mathbf{N}$; \textit{i.e.}, $\varphi $ is
constant on the canonical projections, too.

Then, for each homomorphism $\mu _{i}f\rho _{j}:$ $_{R}R=R_{i}%
\longrightarrow R_{j}=$ $_{R}R$ there exists $r_{ij}\in R$ such that $\mu
_{i}f\rho _{j}(a)=ar_{ij}$ for any $a\in $ $_{R}R=R_{i}$. Therefore, $\mu
_{i}\varphi (f)\rho _{j}(a)=\varphi (\mu _{i})\varphi (f)\varphi (\rho
_{j})(a)=\varphi (\mu _{i}f\rho _{j})(a)=\varphi _{R}(\mu _{i}f\rho
_{j})(a)=a(r_{ij})^{\sigma }=(a^{^{\prime }})^{\sigma }(r_{ij})^{\sigma
}=(a^{^{\prime }}r_{ij})^{\sigma }=(\mu _{i}f\rho _{j}(a^{^{\prime
}}))^{\sigma }$ for any $a\in $ $_{R}R=R_{i}$ and $a^{^{\prime }}=a^{\sigma
^{-1}}$. From this one can easily see that the assignments $\oplus
_{i=1}^{n}R_{i}$ $\ni \oplus _{i=1}^{n}a_{i}\overset{\sigma _{n}^{-1}}{%
\longmapsto }\oplus _{i=1}^{n}a_{i}^{^{\prime }}\overset{f}{\longmapsto }%
f(\oplus _{i=1}^{n}a_{i}^{^{\prime }})\overset{\sigma _{m}}{\longmapsto }%
\sigma _{m}(f(\oplus _{i=1}^{n}a_{i}^{^{\prime }}))\in \oplus _{j=1}^{m}R_{j}
$ define a skew-inner functor $\varphi ^{\sigma }$, and $\varphi ^{\sigma
}=\varphi $.\textit{\ \ \ \ \ \ }$_{\square }\medskip $

\noindent \textbf{Corollary 3.4.} \textit{For any} $\varphi \in Aut$ $(_{R}%
\mathcal{M}_{Sk}^{0})$\textit{\ there exists a skew-inner automorphism }$%
\varphi _{0}\in Aut$ $(_{R}\mathcal{M}_{Sk}^{0})$ \textit{such that} $%
\varphi $ \textit{and} $\varphi _{0}$ \textit{are naturally isomorphic, i.e.,%
} $\varphi \cong $ $\varphi _{0}$ \textit{in the category} $\mathcal{F}(_{R}%
\mathcal{M}_{Sk}^{0},_{R}\mathcal{M}_{Sk}^{0})$ \textit{of endofunctors on} $%
_{R}\mathcal{M}_{Sk}^{0}$\medskip .\textbf{\ }

\noindent \textbf{Proof}. Indeed, as $\varphi \in Aut$ $(_{R}\mathcal{M}%
_{Sk}^{0})$ preserves biproducts, for any $n\in $ $\mathbf{N}$ there are two
mutually inverse isomorphisms $\sum_{i}\varphi (\mu _{i}):\oplus
_{i=1}^{n}R_{i}\longrightarrow \oplus _{i=1}^{n}R_{i}$ and $\sum_{i}\mu
_{i}:\oplus _{i=1}^{n}R_{i}\longrightarrow \oplus _{i=1}^{n}R_{i}$ such that 
$\mu _{i}\sum_{i}\varphi (\mu _{i})=\varphi (\mu _{i})$ and $\varphi (\mu
_{i})\sum_{i}\mu _{i}=\mu _{i}$ for all $i=1,\ldots ,n$. Then, defining $%
\varphi _{0}$ as $\varphi _{0}(f)\overset{def}{=}\sum_{i}\varphi (\mu _{i})$ 
$\varphi (f)$ $\sum_{j}\lambda _{j}$ for any $f:\oplus
_{i=1}^{n}R_{i}\longrightarrow \oplus _{j=1}^{m}R_{j}$ and biproducts $%
\oplus _{i=1}^{n}R_{i}\overset{\pi _{i}}{\underset{\mu _{i}}{%
\rightleftarrows }}R_{i}$ and $\oplus _{j=1}^{m}R_{j}\overset{\rho _{j}}{%
\underset{\lambda _{j}}{\rightleftarrows }}R_{j}$, one may readily verify
that $\varphi _{0}\in Aut$ $(_{R}\mathcal{M}_{Sk}^{0})$\ and always $\varphi
_{0}(\mu _{i})=\mu _{i}$ for all injections $\mu _{i}$. Finally, it is clear
that the isomorphisms $\sum_{i}\varphi (\mu _{i}):\oplus
_{i=1}^{n}R_{i}\longrightarrow \oplus _{i=1}^{n}R_{i}$, $n\in $ $\mathbf{N}$%
, define a natural isomorphism $\varphi _{0}\longrightarrow \varphi $ in $%
\mathcal{F}(_{R}\mathcal{M}_{Sk}^{0},_{R}\mathcal{M}_{Sk}^{0})$.\textit{\ \
\ \ \ \ }$_{\square }\medskip $

It is obvious that skew-inner automorphisms are examples of the more general
notion that we introduce now.\medskip

\noindent \textbf{Definition 3.5.} An automorphism $\varphi \in Aut$ $(_{R}%
\mathcal{M}^{0})$\textbf{\ (}$\varphi \in Aut$ $(_{R}\mathcal{M}_{Sk}^{0})$)
is called \textit{semi-inner} if there exist a semiring automorphism $\sigma
:R\longrightarrow R$ and a family $\{s_{F_{X}}$ $|$ $F_{X}\in |_{R}\mathcal{M%
}^{0}|$ $\}$ ($\{$ $s_{F_{n}}$ $|$ $F_{n}\in |_{R}\mathcal{M}_{Sk}^{0}|\}$)
of monoid isomorphisms $s_{F_{X}}:F_{X}\longrightarrow \varphi (F_{X})$ ( $%
s_{F_{n}}:F_{n}=\oplus _{i=1}^{n}R_{i}\longrightarrow \oplus
_{i=1}^{n}R_{i}=F_{n}$), holding $s_{F_{X}}(ra)=$ $r^{\sigma }s_{F_{X}}(a)$ (%
$s_{F_{n}}(ra)=$ $r^{\sigma }s_{F_{n}}(a)$) for all $r\in R$ and $a\in F_{X}$
($a\in F_{n}$), such that for any $f:F_{X}\longrightarrow F_{Y}$ ($%
f:F_{n}\longrightarrow F_{m}$) the diagrams 
\begin{equation*}
\begin{tabular}{lll}
$F_{X}$ & $\overset{s_{F_{X}}}{\longrightarrow }$ & $\varphi (F_{X})$ \\ 
$f\downarrow $ &  & $\downarrow \varphi (f)$ \\ 
$F_{Y}$ & $\underset{s_{F_{Y}}}{\longrightarrow }$ & $\varphi (F_{Y})$%
\end{tabular}%
\ \text{ \ (}%
\begin{array}{ccc}
F_{n}=\oplus _{i=1}^{n}R_{i} & \overset{s_{F_{n}}}{\longrightarrow } & 
\oplus _{i=1}^{n}R_{i}=F_{n} \\ 
f\downarrow &  & \downarrow \varphi (f) \\ 
F_{m}=\oplus _{j=1}^{m}R_{j} & \underset{s_{F_{m}}}{\longrightarrow } & 
\oplus _{j=1}^{m}R_{j}=F_{m}%
\end{array}%
\text{) \ \ (2)}
\end{equation*}

\noindent are commutative.\medskip

From Corollary 3.4 we immediately obtain\medskip

\noindent \textbf{Corollary 3.6.} \textit{All} $\varphi \in Aut$ $(_{R}%
\mathcal{M}_{Sk}^{0})$ \textit{are semi-inner.\medskip }

\noindent \textbf{Proof}. Indeed, by Corollary 3.4 an automorphism $\varphi
\in Aut$ $(_{R}\mathcal{M}_{Sk}^{0})$ is isomorphic to a skew-inner
automorphism $\varphi _{0}$,$\ $defined, let us say, by a semiring
automorphism $\sigma :R\longrightarrow R$, \textit{i.e.}, $\varphi
_{0}=\varphi _{0}^{\sigma }$. Then, using the notations introduced above,
one can easily verify that the monoid automorphisms $F_{n}=\oplus
_{i=1}^{n}R_{i}\overset{\sigma _{n}}{\longrightarrow }F_{n}=\oplus
_{i=1}^{n}R_{i}\overset{\sum_{i}\varphi (\mu _{i})}{\longrightarrow }%
F_{n}=\oplus _{i=1}^{n}R_{i}$, $F_{n}\in |_{R}\mathcal{M}_{Sk}^{0}|$, $n\in $
$\mathbf{N}$, define the needed monoid automorphisms $s_{F_{n}}$, $F_{n}\in
|_{R}\mathcal{M}_{Sk}^{0}|$, $n\in $ $\mathbf{N}$, in Definition 3.5.\textit{%
\ \ \ \ \ \ }$_{\square }\medskip $

By Lemma 2.5, any automorphism $\overline{\varphi }\in Aut$ $(_{R}\mathcal{M}%
_{Sk}^{0})$ induces the automorphism $\overline{\varphi }^{i}\in Aut$ $(_{R}%
\mathcal{M}^{0})$, about which the following observation is true.\medskip

\noindent \textbf{Lemma 3.7.} \textit{For any} $\overline{\varphi }\in Aut$ $%
(_{R}\mathcal{M}_{Sk}^{0})$\textit{, the automorphism} $\overline{\varphi }%
^{i}\in Aut$ $(_{R}\mathcal{M}^{0})$ \textit{is semi-inner.}\medskip

\noindent \textbf{Proof}. By Corollary 3.6, an automorphism $\overline{%
\varphi }\in Aut$ $(_{R}\mathcal{M}_{Sk}^{0})$ is semi-inner. So, let $%
\{s_{F_{n}}:F_{n}=\oplus _{i=1}^{n}R_{i}\longrightarrow \oplus
_{i=1}^{n}R_{i}=F_{n},F_{n}\in |_{R}\mathcal{M}_{Sk}^{0}|$, $n\in $ $\mathbf{%
N}\}$ be the corresponding monoid automorphisms from Definition 3.5. Also,
as $_{R}\mathcal{M}_{Sk}^{0}$ is a skeleton of $_{R}\mathcal{M}^{0}$, there
is a family $\{i_{F_{X}}:F_{X}\longrightarrow F_{n}|$ $F_{X}\in |_{R}%
\mathcal{M}^{0}|,|X|=n$, and $i_{F_{n}}=1_{F_{n}}$ if $F_{X}=F_{n}\in |_{R}%
\mathcal{M}_{Sk}^{0}|\}$ of semimodule isomorphisms. From this it
immediately follows that $\{s_{F_{X}}:F_{X}\overset{i_{F_{X}}}{%
\longrightarrow }F_{n}\overset{s_{F_{n}}}{\longrightarrow }F_{n}\overset{%
i_{F_{X}}^{-1}}{\longrightarrow }F_{X}$ $|$ $F_{X}\in |_{R}\mathcal{M}%
^{0}|\} $ is a family of monoid automorphisms such as is requested in
Definition 3.5.\textit{\ \ \ \ \ \ }$_{\square }\medskip $

\noindent \textbf{Lemma 3.8.}\textit{\ Let }$\varphi ,\psi \in Aut$ $(_{R}%
\mathcal{M}^{0})$ \textit{be semi-inner. Then their composite }$\varphi \psi
\in Aut$ $(_{R}\mathcal{M}^{0})$ \textit{is also semi-inner.\medskip }

\noindent \textbf{Proof}. Let $\{s_{F_{X}}:F_{X}\longrightarrow \varphi
(F_{X})$ $|$ $F_{X}\in |_{R}\mathcal{M}^{0}|$ $\}$ and $\sigma
:R\longrightarrow R$ and $\{t_{F_{X}}:F_{X}\longrightarrow \varphi (F_{X})$ $%
|$ $F_{X}\in |_{R}\mathcal{M}^{0}|\}$ and $\tau :R\longrightarrow R$
correspond to $\varphi $ and $\psi $, respectively. Then, it can be easily
verified that the family $\{F_{X}\overset{s_{F_{X}}}{\longrightarrow }%
\varphi (F_{X})\overset{t_{\varphi (F_{X})}}{\longrightarrow }\psi (\varphi
(F_{X}))$ $|$ $F_{X}\in |_{R}\mathcal{M}^{0}|\}$ of monoid isomorphisms and $%
\sigma \tau :R\longrightarrow R$ set $\varphi \psi $ to be semi-inner, too.%
\textit{\ \ \ \ \ \ }$_{\square }\medskip $

\noindent \textbf{Theorem 3.9.} \textit{Let} $R$ \textit{be an IBN-semiring.
Then all automorphisms }$\varphi \in Aut$ $(_{R}\mathcal{M}^{0})$ \textit{%
are semi-inner.\medskip }

\noindent \textbf{Proof}. Let $\varphi \in Aut$ $(_{R}\mathcal{M}^{0})$. By
Proposition 2.9, $\varphi \in EqnAut$ $(_{R}\mathcal{M}^{0})=Aut$ $(_{R}%
\mathcal{M}^{0})$, and, therefore, by Proposition 2.7, $\varphi =\varphi
_{S}\varphi _{I}$, where $\varphi _{S}\in StAut$\textit{\ }$(_{R}\mathcal{M}%
^{0})$ and $\varphi _{I}\in Int$ $(_{R}\mathcal{M}^{0})$.

By Proposition 2.6 and Lemma 3.7, $\varphi _{S}$ is naturally isomorphic to
the semi-inner automorphism $\overline{\varphi _{S}}^{i}\in Aut$ $(_{R}%
\mathcal{M}^{0})$. Then, if $\{s_{F_{X}}:F_{X}\longrightarrow F_{X}$ $|$ $%
F_{X}\in |_{R}\mathcal{M}^{0}|\}$ and $\sigma :R\longrightarrow R$ define
the semi-inner \ $\overline{\varphi _{S}}^{i}$, and a family $\{\eta
_{F_{X}}:F_{X}\longrightarrow F_{X}$ $|$ $F_{X}\in |_{R}\mathcal{M}^{0}|\}$
of semimodule automorphisms define a natural isomorphism $\overline{\varphi
_{S}}^{i}\overset{\bullet }{\longrightarrow }$ $\varphi _{S}$, it is readily
verified that $\{F_{X}\overset{s_{F_{X}}}{\longrightarrow }F_{X}\overset{%
\eta _{F_{X}}}{\longrightarrow }F_{X}$ $|$ $F_{X}\in |_{R}\mathcal{M}^{0}|\}$
and $\sigma :R\longrightarrow R$ make $\varphi _{S}$ semi-inner, too. From
this, the obvious fact that inner automorphisms are semi-inner, and Lemma
3.8 we conclude our proof.\textit{\ \ \ \ \ \ }$_{\square }\medskip $

As an immediate corollary of Theorem 3.9 we obtain the following
result.\medskip

\noindent \textbf{Theorem 3.10.} \textit{Let }$R$\textit{\ be a semiring of
one of the following classes of semirings: 1) Artinian (left or right)
rings; 2) Noetherian (left or right) rings; 3) Commutative rings; 4)
PI-rings; 5) Additively-idempotent division semirings, in particular,
schedule algebras; \ 6) Division semirings. Then all automorphisms }$\varphi
\in Aut$ $(_{R}\mathcal{M}^{0})$ \textit{are semi-inner.\medskip }

\noindent \textbf{Proof}. All rings of the first four classes have IBN (see,
for example,\ \cite{lam:lomar} and/or \cite{rowen:rt}). Semirings of the
last two classes also are an IBN-semirings by \cite[Theorem 5.3]{hebwei:ros}.%
\textit{\ \ \ \ \ \ }$_{\square }\medskip $

Following \cite[Definition 3.2]{mashplts:aocfa} and taking a variety $\Theta 
$ to be \textit{perfect} provided $Aut$ $(\Theta ^{0})=Int$ $(\Theta ^{0})$,
we have the following corollary of Theorem 3.9. \medskip

\noindent \textbf{Corollary 3.11.} \textit{If the group of automorphisms of
an IBN-semiring} $R$ \textit{is trivial, then the variety} $_{R}\mathcal{M}$ 
\textit{is perfect.\ \ \ \ \ \ }$_{\square }\medskip $

Then, from this result we obtain \medskip

\noindent \textbf{Corollary 3.12.} \textit{The varieties of abelian groups,
abelian monoids, and commutative Clifford monoids are perfect.\medskip }

\noindent \textbf{Proof}. Concerning the variety of abelian groups, the
statement is obvious.

As to the variety of abelian monoids, the result follows from the
observation that this variety coincides with the variety of semimodules over
the semiring $\mathbf{N}_{0}$, which has IBN by \cite[Theorem 5.3]%
{hebwei:ros}.

Regarding the variety of commutative Clifford monoids, from %
\cite[Propositions 10 and 11]{katsov:tpies}$\ $\ and \cite[Theorem 5.3]%
{hebwei:ros}, we have that this variety also coincides with a variety of
semimodules over an IBN-semiring with a trivial automorphism group.\textit{\
\ \ \ \ \ }$_{\square }\medskip $

Now, if a variety $\Theta $ is a category $_{R}\mathcal{M}$ of left
semimodules over a semiring $R$, by Proposition 2.10, any automorphism $%
\varphi :$ $_{R}\mathcal{M}^{0}\longrightarrow $ $_{R}\mathcal{M}^{0}$%
\textit{\ }is equinumerous provided that in the corresponding monogenic
semigroup $M(n,h)$ the index $n>1$, or $n=1$ and $h$ is\textit{\ }either%
\textit{\ }$1$ or $2$.\ From this observation and following closely the
scheme of the proof of Theorem 3.9, one can readily extend that
result.\medskip

\noindent \textbf{Theorem 3.13.}\ \textit{Let} $R$ \textit{be an
IBN-semiring, or }$n>1$\textit{, or} $n=1$ \textit{and} $h$ \textit{be
either }$1$\textit{\ or} $2$\textit{\ for a monogenic semigroup} $M(n,h)$%
\textit{. Then all automorphisms }$\varphi \in Aut$ $(_{R}\mathcal{M}^{0})$ 
\textit{are semi-inner.\medskip\ \ \ \ \ \ }$_{\square }$

Analyzing the proofs of Theorems 3.9 and 3.13, one may easily make the
following observation.\medskip

\noindent \textbf{Proposition 3.14.} \textit{If for a semiring} $R$ \textit{%
all} \textit{automorphisms }$\varphi \in Aut$ $(_{R}\mathcal{M}^{0})$ 
\textit{are equinumerous, then they are semi-inner.\ \ \ \ \ \ }$_{\square
}\medskip $

In light of Proposition 3.14, Theorem 3.13 and Problem 1, we have the
following interesting open problems.\medskip

\noindent \textbf{Problem 2.} Does there exist a (semi)ring $R$ of a type $%
(1,h)$, $h>2$, having a non-semi-inner automorphism $\varphi \in Aut$ $(_{R}%
\mathcal{M}^{0})$? (Our conjecture is that the answer is ``Yes.'')\medskip

\noindent \textbf{Problem 3.} Does there exist a (semi)ring $R$ of a type $%
(1,h)$, $h>2$, having a semi-inner non-equinumerous automorphism $\varphi
\in Aut$ $(_{R}\mathcal{M}^{0})$? (Our conjecture is that the answer is
``Yes.'')\medskip

\noindent \textbf{Problem 4. }Is the converse of Proposition 3.14 true? \
(Our conjecture is that the answer is ``No.'')\medskip

Finally, the following interesting observations in respect to a group $Out$ $%
(_{R}\mathcal{M}^{0})=$ $Aut$ $(_{R}\mathcal{M}^{0})/Int$ $(_{R}\mathcal{M}%
^{0})$ of outer automorphisms of the category $_{R}\mathcal{M}^{0}$ will
conclude this section. \medskip

\noindent \textbf{Theorem 3.15.} \textit{Let }$Out$ $(R)=$ $Aut$ $(R)/Int$ $%
(R)$ \textit{denote the group of outer automorphisms of a semiring} $R$%
\textit{, and all automorphisms of} $Aut$ $(_{R}\mathcal{M}^{0})$ \textit{be
equinumerous. Then,} $Out$ $(_{R}\mathcal{M}^{0})\cong Out$ $(R)$\textit{%
.\medskip }

\noindent \textbf{Proof}. By Proposition 2.7, $Aut$ $(_{R}\mathcal{M}%
^{0})=StAut$\textit{\ }$(_{R}\mathcal{M}^{0})Int$ $(_{R}\mathcal{M}^{0})$,
and, hence, $Out$ $(_{R}\mathcal{M}^{0})=$ $Aut$ $(_{R}\mathcal{M}^{0})/Int$ 
$(_{R}\mathcal{M}^{0})\cong $ $StAut$\textit{\ }$(_{R}\mathcal{M}^{0})Int$ $%
(_{R}\mathcal{M}^{0})/Int$ $(_{R}\mathcal{M}^{0})\cong $ $StAut$\textit{\ }$%
(_{R}\mathcal{M}^{0})/(StAut\mathit{\ }(_{R}\mathcal{M}^{0})\cap Int$ $(_{R}%
\mathcal{M}^{0}))$. Then, using the same arguments and notations as in
Proposition 3.3, one can easily see that the assignment $StAut$\textit{\ }$%
(_{R}\mathcal{M}^{0})\ni \varphi \longmapsto \varphi _{R}\longmapsto \sigma
\in Aut$ $(R)$ defines a group isomorphism

\noindent $StAut$\textit{\ }$(_{R}\mathcal{M}^{0})\longrightarrow Aut$ $(R)$%
; and moreover, if $\varphi \in StAut\mathit{\ }(_{R}\mathcal{M}^{0})\cap
Int $ $(_{R}\mathcal{M}^{0})$, then there is a family $\{i_{F_{X}}:F_{X}%
\longrightarrow F_{X}$ $|$ $F_{X}\in |_{R}\mathcal{M}^{0}|\}$ of semimodule
isomorphisms such that $\varphi (f)=i_{F_{Y}\text{ }}fi_{F_{X}}^{-1}$ for
any $f:F_{X}\longrightarrow F_{Y}$, and, therefore, $\varphi \longmapsto
\sigma =i_{R}\in Aut$ $(R)$ under the isomorphism $StAut$\textit{\ }$(_{R}%
\mathcal{M}^{0})\longrightarrow Aut$ $(R)$. As a result of that, we have
that $Out$ $(_{R}\mathcal{M}^{0})$ $\cong $ $StAut$\textit{\ }$(_{R}\mathcal{%
M}^{0})/(StAut\mathit{\ }(_{R}\mathcal{M}^{0})\cap Int$ $(_{R}\mathcal{M}%
^{0}))\cong Aut$ $(R)/Int$ $(R)=Out$ $(R)$.\textit{\ \ \ \ \ \ }$_{\square
}\medskip $

\noindent \textbf{Corollary 3.16.}\textit{\ Let }$R$ \textit{be an
IBN-semiring, in particular, a semiring of one of the classes of semirings
of Theorem 3.10. Then the group of outer automorphisms of the category} $_{R}%
\mathcal{M}^{0}$ \textit{is isomorphic to the group of outer automorphisms
of the semiring }$R$\textit{.\ \ \ \ \ \ }$_{\square }\medskip $

\section{Automorphisms of Categories of Free Lie Modules and Free Restricted
Lie Modules}

From now on, let $K$ be a commutative associative ring with unity $1$, a Lie
algebra $L$ over $K$ a free $K$-module with a totally ordered $K$-basis $%
E=\{e_{j}|j\in J\}$, and $U(L)$ its universal enveloping algebra. And thus,
by the Poincare-Birkhoff-Witt Theorem (PBWT) (see, for example, %
\cite[Theorem 2.5.3]{bah:irila}, or \cite[Theorem 3.4.3]{serre:ser}), we may
consider $L$ to be a subalgebra of the commutator algebra $U^{(-)}$ of the
algebra $U(L)$.

Let $_{L}\mathcal{M}$ and $_{U(L)}\mathcal{M}$ denote the categories of $L$-
and $U(L)$-modules, respectively. It is well known (see, \textit{e.g.}, %
\cite[Section 2.5.5]{bah:irila}) that these categories are isomorphic in
such a way that free $L$-modules generated by a set of free generators $%
X=\{x_{1},...,x_{n},....\}$ go to free $U(L)$-modules generated by the same $%
X$. Using this fact, we first consider automorphisms of the full subcategory 
$_{U(L)}\mathcal{M}^{0}$ of finitely generated free $U(L)$-modules $%
F_{X},\;X\subseteq \mathcal{U},\;|X|<\infty $, of the category $_{U(L)}%
\mathcal{M}$; and then, we obtain a description of automorphisms of the full
subcategory $_{L}\mathcal{M}^{0}$ of finitely generated free $L$-modules $%
F_{X},\;X\subseteq \mathcal{U},\;|X|<\infty $, of the category $_{L}\mathcal{%
M}$.\medskip\ 

We will use the following description of a basis over $K$ ($K$-basis) of a
free $L$-module $F_{X}$ $\in |_{L}\mathcal{M}|$.\medskip

\noindent \textbf{Theorem 4.1.}(cf. \cite[Theorem 1.6.11.]{bah:irila}) 
\textit{\ Let} $F_{X}$ \textit{be a free }$L$-\textit{module with free
generators} $X=\{x_{1},...,x_{n},....\}$\textit{,} \textit{and} $%
E=\{e_{1},...,e_{n},...\}$ \textit{a totally ordered }$K$\textit{-basis of} $%
L$\textit{.} \textit{Then the set }$\{x_{1},...,x_{n},....\}\cup
\{e_{1}e_{2}...e_{n}x_{i}|$ $x_{i}\in X,\;e_{1},e_{2},...,e_{n}\in E,$ $%
e_{1}\geq e_{2}\geq ...\geq e_{n},\;n=1,2,...\}$ \textit{is a }$K$\textit{%
-basis of }$F_{X}$. \medskip

Using this description, we have the following result.\medskip

\noindent \textbf{Proposition 4.2.} \textit{Let a Lie algebra }$L$ \textit{%
be a free }$K$\textit{-module with a totally ordered }$K$\textit{-basis }$%
E=\{e_{1},...,e_{n},...\}$. \textit{Then} $_{L}\mathcal{M}$ \textit{and} $%
_{U(L)}\mathcal{M}$ \textit{are IBN-varieties.} \medskip

\noindent \textbf{Proof.} Because of the isomorphism of the categories $_{L}%
\mathcal{M}$ and $_{U(L)}\mathcal{M}$, it is enough to show that $_{L}%
\mathcal{M}$ and is an\textit{\ }IBN-variety. So, let $F_{X}$ and $F_{Y}$ be
two isomorphic free $L$-modules generated by $X=\{x_{1},...,x_{n}\}$ and $%
Y=\{y_{1},...,y_{m}\}$, respectively, and $\varphi :F_{X}\rightarrow F_{Y}\ $%
an $L$-module isomorphism between $F_{X}$ and $F_{Y}$. Also, by Theorem 4.1
there are the submodules $M_{1}\subseteq $ $F_{X}$ with the $K$-basis $%
S_{1}=\{e_{1}e_{2}...e_{s_{1}}x_{i}|$ $s_{1}\geq 1,\;x_{i}\in
X,\;e_{1},e_{2},...,e_{s_{1}}\in E\}$, and $M_{2}\subseteq $ $F_{Y}$ with
the $K$-basis $S_{2}=\{e_{1}e_{2}...e_{s_{2}}y_{i}|$ $s_{2}\geq 1,\;y_{i}\in
Y,\;e_{1},e_{2},...,e_{s_{2}}\in E\}$.

It is obvious that $\varphi (M_{1})\subseteq M_{2}$, and the $L$-modules $%
F_{X}/M_{1}$ and $F_{Y}/\varphi (M_{1})$ are isomorphic. In fact, $\varphi
(M_{1})=M_{2}$. Indeed, if $\varphi (M_{1})\subset M_{2}$, \textit{i.e.}, $%
\varphi (M_{1})$ were a proper subset of $M_{2}$, then it is easy to see
that the algebra $L$ would act on $F_{Y}/\varphi (M_{1})$ non-trivially,
however, its actions on $F_{X}/M_{1}$ are trivial.

Thus, $F_{X}/M_{1}$ and $F_{Y}/M_{2}$ are isomorphic free $K$-modules having
the finite bases $X$ and $Y$, respectively, and, since a commutative ring $K$
is an IBN-ring, we have $n=m$.\textit{\ \ \ \ \ \ }$_{\square }\medskip $

\noindent \textbf{Proposition 4.3.} \textit{Let} $L$ \textit{be a Lie
algebra over an integral domain} $K$. \textit{Then all units of} $U(L)$ 
\textit{belong to} $K$. \medskip

\noindent \textbf{Proof.} We give a proof that is quite similar to that that
the algebra $U(L)$ over a field $F$ is a domain (see, for example, %
\cite[Theorem 5.6.]{jacobson:jac}).

The graded algebra $gr(U)$ associated with $U=U(L)$ can be described as
follows. Let $L$ be a Lie algebra and $E=\{e_{1},...,e_{n},...\}$ its
totally ordered basis. Define submodules $U_{n},\;n\in N$ of $U(L)$ as $%
U_{n}=K+L+LL+...+L^{n}$, where $L^{i}$ is the submodule of $U(L)$ generated
by all products of $i$ elements of the algebra $L$. For each $n\in N$ define 
$K$-module $U^{n}$ by $U^{n}=U_{n}/U_{n-1}$, where by a convention $U_{-1}=0$%
. Then, $gr(U)=\oplus _{n\in N}U^{n}$ is the graded algebra associated with
the filtration $U_{n}$ of $U(L)$. It is well known (see, \textit{e.g.}, %
\cite[Proposition 3.4.1 and Theorem 3.4.2]{serre:ser}) that $gr(U)\simeq K[{%
e_{j}}]$, where $K[{e_{j}}]$ is the polynomial algebra in the commuting
variables $e_{j}\in E$. Since a polynomial algebra is a domain, it follows
that $gr(U)$ also is a domain.

Let $u$ be a unit in $U(L)$, and $u\cdot u^{-1}=1$ for $u^{-1}\in U(L)$. If $%
u\notin U_{0}=K$, then $u^{-1}\notin U_{0}$, too. Hence, there exist two
unique natural numbers $n\geq 1$ and $m\geq 1$ such that $u\in U_{n}$, $%
u\notin U_{n-1}$, and $u^{-1}\in U_{m}$, $u^{-1}\notin U_{m-1}$. Now
consider $\bar{u}=u+U_{n-1}$ and $\overline{u^{-1}}=u^{-1}+U_{m-1}$ in $%
gr(U) $. They are nonzero homogeneous elements of $gr(U)$ of degree $n$ and $%
m$, respectively. Since $gr(U)$ is a domain, $\bar{u}\cdot \overline{u^{-1}}$
is a nonzero element of degree $n+m$. However, $\bar{u}\cdot \overline{u^{-1}%
}=1+U_{n+m-1}=U_{n+m-1}$. This contradiction proves that $u\in U_{0}=K$. 
\textit{\ \ \ \ \ \ }$_{\square }\medskip $

\noindent \textbf{Corollary 4.4.} \textit{Let} $L$ \textit{be a Lie algebra
as above,} \textit{and} $\varphi :U(L)\rightarrow U(L)$ \textit{a ring
automorphism of the universal enveloping algebra} $U(L)$. \textit{Then,} $%
\varphi (K)=K$. \medskip

\noindent \textbf{Proof.} Let $\tilde{K}=Frac\,K$ be the field of fractions
of $K$. Denote by $L_{\tilde{K}}$ and $U(L)_{\tilde{K}}$ the algebras
derived from $L$ and $U(L)$ by the extension of $K$. Let $\alpha \in \tilde{K%
},\;\alpha \neq 0$. Then, $\varphi (\alpha )$ is a unit of $U(L)_{\tilde{K}}$%
. By Proposition 4.3, $\varphi (\alpha )\in \tilde{K}$. Hence, $\varphi (%
\tilde{K})=\tilde{K}$, and thus, $\varphi (K)=K$. \textit{\ \ \ \ \ \ }$%
_{\square }\medskip $

\noindent \textbf{Definition 4.5.} For (associative, or Lie) algebras\textit{%
\ }$A_{1}$ and\textit{\ }$A_{2}$\textit{\ }over $K$, an automorphism $\delta 
$ of $K$, and a ring homomorphism $\varphi :A_{1}\rightarrow A_{2}$, a pair $%
\phi =(\delta ,\varphi )$ is called a \textit{semi-morphism} from $A_{1}$ to 
$A_{2}$ if $\varphi (\alpha \cdot u)=\alpha ^{\delta }\cdot \varphi (u)$ for
any $\alpha \in K,\;u\in A_{1}$\medskip

\noindent \textbf{Definition 4.6.}\textit{\ }An automorphism $\varphi \in
Aut(_{U(L)}M^{0})$ is called \textit{semi-inner} if there exist a
semi-automorphism $\chi =(\tau ,\tilde{\varphi}):U(L)\longrightarrow U(L)$,
where $\tau \in Aut\;K$ and $\tilde{\varphi}:U(L)\rightarrow U(L)$ is a ring
automorphism of $U(L)$, and a family $\{s_{F_{X}}|$ $s_{F_{X}}:F_{X}%
\rightarrow \varphi (F_{X}),\;F_{X}\in |_{U(L)}M|\}$ of monoid isomorphisms
holding $s_{F_{X}}(\alpha lu)=\alpha ^{\tau }l^{\tilde{\varphi}}s_{F_{X}}(u)$
for all $l\in U(L),\;\alpha \in K$ and $u\in F_{X}$, such that for any $%
f:F_{X}\longrightarrow F_{Y}$ the diagram 
\begin{equation*}
\begin{tabular}{lll}
$F_{X}$ & $\overset{s_{F_{X}}}{\longrightarrow }$ & $\varphi (F_{X})$ \\ 
$f\downarrow $ &  & $\downarrow \varphi (f)$ \\ 
$F_{Y}$ & $\underset{s_{F_{Y}}}{\longrightarrow }$ & $\varphi (F_{Y})$%
\end{tabular}%
\end{equation*}%
\noindent is commutative.\medskip

\noindent \textbf{Theorem 4.7.} \textit{Let } $L$ \textit{be a Lie algebra
over an integral domain} $K$. \textit{Then all automorphisms } $\varphi \in
Aut(_{U(L)}\mathcal{M}^{0})$ \textit{are semi-inner.}\medskip

\noindent \textbf{Proof.} Let us consider an algebra $U(L)$ as a ring, and
denote by $_{RU(L)}\mathcal{M}^{0}$ the category of finitely generated free
modules over the ring $U(L)$. By proposition 4.2, $U(L)$ is an IBN-algebra,
and, hence, it is easy to see that $U(L)$ is also an IBN-ring. Therefore, by
Theorem 3.9 all automorphisms of $_{RU(L)}\mathcal{M}^{0}$ are semi-inner.
And since by Corollary 4.4 $\sigma (K)=K$ for any ring automorphism $\sigma
:U(L)\rightarrow U(L)$, all automorphisms of category $_{U(L)}\mathcal{M}%
^{0} $ are semi-inner, too. \textit{\ \ \ \ \ \ }$_{\square }\medskip $

Obviously modifying the notion of a semi-inner automorphism for the category 
$_{L}\mathcal{M}^{0}$, and taking into account the isomorphism of the
categories $_{L}\mathcal{M}$ and $_{U(L)}\mathcal{M}$, from Theorem 4.7 we
obtain \medskip

\noindent \textbf{Corollary 4.8.} \textit{Let } $L$ \textit{be a Lie algebra
over an integral domain} $K$. \textit{Then, }$Aut(_{L}\mathcal{M}%
^{0})=Aut(_{U(L)}\mathcal{M}^{0})$\textit{, and} \textit{all automorphisms} $%
\varphi \in Aut(_{L}\mathcal{M}^{0})$ \textit{are semi-inner.} \textit{\ \ \
\ \ \ }$_{\square }\medskip $

Now let $p$ be a prime number, $\Phi $ an associative commutative ring with
unity $1$, and $p\cdot \Phi =0$. Recall \cite[Section 1.11]{bah:irila} that
a restricted Lie algebra (or $p$-Lie algebra) is a Lie algebra $G$ over $%
\Phi $ with an unary operation $g\longmapsto g^{p}$,$\;g\in G$, satisfying
the following identities:

(1) $(\lambda g)^{p}=\lambda^{p}g^{p} $, for all $\lambda \in \Phi,\; g\in G$%
;

(2) $ad\ g^{p}=(ad\ g\ )^{p}$, $g\in G$;

(3) $(g_{1}+g_{2})^{p}=g_{1}^{p}+g_{2}^{p}+\Sigma
_{i=1}^{p-1}s_{i}(g_{1},g_{2}),\;g_{1},g_{2}\in G$, where $%
s_{i}(g_{1},g_{2}) $ is the coefficient of $t^{i-1}$ in $ad$ $%
(tg_{1}+g_{2})^{p-1}(g_{1})$.

For a restricted Lie algebra $G$, the notions of an enveloping algebra and
of an universal enveloping algebra $U_{p}(G)$ are defined in the same
fashion as their corresponding Lie algebra analogs. Also, if $G$ is a free $%
\Phi $-module, and $U_{p}(G)$ is its universal enveloping algebra, then the
PBWT is valid for a $p$-Lie algebra $G$, too. In particular, the $p$-Lie
algebra $G$ is a subalgebra of the commutator algebra $U_{p}^{(-)}(G)$ of $%
U_{p}(G)$.

Now considering the categories $_{G}\mathcal{M}_{p}^{0}$ and $_{U_{p}(G)}%
\mathcal{M}_{p}^{0}$ of free $p$-Lie modules over $G$ and free modules over $%
U_{p}(G)$, respectively, one can easily obtain the following analogue of
Theorem 4.1. \medskip

\noindent \textbf{Theorem 4.9.} \textit{\ Let }$F_{X}$\textit{\ be a free }$%
p $\textit{-Lie module over a }$p$\textit{-Lie algebra }$G$\textit{\ with
free generators }$X=\{x_{1},...,x_{n},....\}$\textit{, and }$%
E=\{e_{1},...,e_{n},...\}$\textit{\ a totally ordered }$\Phi $\textit{-basis
of }$G$\textit{. Then the set }$\{x_{1},...,x_{n},....\}$ $\cup $

\noindent $\{e_{1}^{\alpha _{1}}e_{2}^{\alpha _{2}}...e_{n}^{\alpha
_{n}}x_{i}|$\textit{\ }$x_{i}\in X,\;e_{1},e_{2},...,e_{n}\in E,$\textit{\ }$%
e_{1}\geq e_{2}\geq ...\geq e_{n},\;n=1,2,...;\;0\leq \alpha _{i}\leq p-1\}$%
\textit{\ is a }$\Phi $\textit{-basis of }$F_{X}$\textit{.} \textit{\ \ \ \
\ \ }$_{\square }\medskip $ \medskip

Then, appropriately modifying the notion of a semi-inner automorphism for
the categories $_{G}\mathcal{M}_{p}^{0}$ and $_{U_{p}(G)}\mathcal{M}_{p}^{0}$
and following the same scheme as in the proof of Theorem 4.7, we establish
\medskip

\noindent \textbf{Theorem 4.10.} \textit{Let }$G$\textit{\ be an p-Lie
algebra over an integral domain }$\Phi $\textit{. Then all automorphisms of }%
$Aut(_{G}M_{p}^{0})$\textit{\ and of }$Aut(_{U_{p}(G)}M_{p}^{0})$\textit{\
are semi-inner. \ \ \ \ \ \ }$_{\square }\medskip $ \medskip

Finally, we wish to consider some connections between the groups $Aut(_{U(L)}%
\mathcal{M}^{0})\ $and $Aut\,U(L)$, and to do that we will need the
following lemma.\medskip

\noindent \textbf{Lemma 4.11.} \textit{Let }$L$\textit{\ be a Lie algebra
over an integral domain }$K$\textit{, and }$F_{1}$\textit{\ a free cyclic }$%
U(L)$\textit{-module with a free generator }$x_{1}$\textit{, i.e., }$%
F_{1}=U(L)x_{1}$\textit{. Then, }$Aut\,F_{1}=Aut\,Kx_{1}$\textit{, where }$%
Aut\,F_{1}$\textit{\ is the group of all }$U(L)$\textit{-module
automorphisms of }$F_{1}$\textit{, and }$Kx_{1}$\textit{\ is a cyclic }$K$%
\textit{-module generated by }$x_{1}$\textit{.\medskip }

\noindent \textbf{Proof.} Let $\theta :U(L)x_{1}\rightarrow U(L)x_{1}$ be a $%
U(L)$-module automorphism of $F_{1}$ such that $\theta (x_{1})=ux_{1},\,u\in
U(L)$. Let $\theta ^{-1}(x_{1})=vx_{1},\,v\in U(L)$. Then $\theta \theta
^{-1}(x_{1})=x_{1}=uvx_{1}$, \textit{i.e.}, $uv=1$. By Proposition 4.3, $%
u,\,v\in K^{\ast }\subseteq K$, where $K^{\ast }$ is the group of all units
of $K$. This completes the proof. \medskip \textit{\ \ \ \ \ \ }$_{\square }$

Let $\varphi \in StAut(_{U(L)}\mathcal{M}^{0})$, and $\varphi _{Ux_{1}}$ be
a restriction of $\varphi $ on $Ux_{1}$. If $\nu _{l}:Ux_{1}\rightarrow
Ux_{1}$ is an endomorphism of $Ux_{1}$ such that $\nu
_{l}(x_{1})=lx_{1},\;l\in U(L)$, then $\varphi _{Ux_{1}}(\nu _{l})=\nu
_{l^{\sigma }}$ for some mapping $\sigma :U(L)\rightarrow U(L)$
corresponding to $\varphi $, and we have the following observation.\medskip

\noindent \textbf{Proposition 4.12.} \textit{If }$\sigma |_{K}=\sigma _{K}$%
\textit{\ for an integral domain }$K$\textit{, then }$\sigma
:U(L)\rightarrow U(L)$\textit{\ is a ring automorphism of }$U(L)$\textit{,
and }$\sigma _{K}\in Aut\,K$\textit{. \medskip }

\noindent \textbf{Proof.} As $_{U(L)}\mathcal{M}^{0}$ is clearly a
semi-additive category with biproducts and a zero object, by %
\cite[Proposition 12.2.5 and Convention 12.2.6]{schu:cat} the functor $%
\varphi \in Aut(_{U(L)}\mathcal{M}^{0})$ is additive. Then, $\varphi
(\upsilon _{l_{1}}+\upsilon _{l_{2}})=\varphi (\upsilon
_{l_{1}+l_{2}})=\upsilon _{(l_{1}+l_{2})^{\sigma }}$ for any elements $%
l_{1},l_{2}\in U(L)$. On the other hand $\varphi (\upsilon _{l_{1}}+\upsilon
_{l_{2}})=\varphi (\upsilon _{l_{1}})+\varphi (\upsilon _{l_{2}})=\nu
_{l_{1}^{\sigma }}+\nu _{l_{2}^{\sigma }}=\nu _{l_{1}^{\sigma
}+l_{2}^{\sigma }}.$ Thus, $(l_{1}+l_{2})^{\sigma }=l_{1}^{\sigma
}+l_{2}^{\sigma }$.

Similarly, we have $\varphi (\nu _{l_{1}}\nu _{l_{2}})(x_{1})=(\varphi (\nu
_{l_{1}})\varphi (\nu _{l_{2}}))(x_{1})$

\noindent $=\varphi (\nu _{l_{1}})(\varphi (\nu _{l_{2}})(x_{1}))=\nu
_{l_{1}^{\sigma }}(\nu _{l_{2}^{\sigma }}(x_{1}))=l_{1}^{\sigma
}(l_{2}^{\sigma }(x_{1}))=\nu _{l_{1}^{\sigma }l_{2}^{\sigma }}(x_{1})$, and 
$\varphi (\nu _{l_{1}}\nu _{l_{2}})(x_{1})=\varphi (\nu
_{l_{1}l_{2}})(x_{1})=\nu _{(l_{1}l_{2})^{\sigma }}(x_{1})$. Thus, $%
(l_{1}l_{2})^{\sigma }=l_{1}^{\sigma }l_{2}^{\sigma }$, and we have proved
that $\sigma $ is a ring homomorphism of $U(L)$. Since $\varphi \in
Aut(_{U(L)}\mathcal{M}^{0})$, the map $\sigma $ is a ring automorphism of $%
U(L)$; and, by Corollary 4.4, $\sigma _{K}(K)=K$, therefore, $\sigma _{K}$
is an automorphism of $K$. \textit{\ \ \ \ \ \ }$_{\square }\medskip $

Let $\sigma _{K}$ and $\sigma $ be the automorphisms of the integral domain $%
K$ and of the ring $U(L)$, respectively. Define a mapping $\sigma
_{X}:F_{X}\rightarrow F_{X},\;F_{X}\in $ $|_{U(L)}\mathcal{M}^{0}|$ in the
following way: if $u=\sum_{i=1}^{n}\alpha _{i}u_{i}x_{i}$, where $u_{i}\in
U(L),\;\alpha _{i}\in K,\;x_{i}\in X$, then $\sigma
_{X}(u)=\sum_{i=1}^{n}(\alpha _{i}u_{i})^{\sigma }x_{i}=$

\noindent $\sum_{i=1}^{n}\alpha _{i}^{\sigma _{K}}u_{i}^{\sigma }x_{i}$. It
is evident that for all $u,v\in F_{X}$, $l\in U(L)$ and $\alpha \in K$, we
have $\sigma _{X}(u+v)=\sigma _{X}(u)+\sigma _{X}(v)$ and $\sigma
_{X}(\alpha lu)=\alpha ^{\sigma _{K}}l^{\sigma }\sigma _{X}(u)$.\medskip

\noindent \textbf{Lemma 4.13.} \textit{Let }$\nu :F_{X}\rightarrow F_{Y}$%
\textit{\ be a }$U(L)$\textit{-module homomorphism, and }$\sigma
_{X}:F_{X}\rightarrow F_{X}$\textit{, }$\sigma _{Y}:F_{Y}\rightarrow F_{Y}$%
\textit{. Then }$\sigma _{Y}\nu \sigma _{X}^{-1}:F_{X}\rightarrow F_{Y}$%
\textit{\ is a }$U(L)$\textit{-module homomorphism, too.\medskip }

\noindent \textbf{Proof.} We have $\sigma _{Y}\nu \,\sigma
_{X}^{-1}(lu)=\sigma _{Y}\nu (l^{\sigma ^{-1}}\sigma _{X}^{-1}(u))=\sigma
_{Y}(l^{\sigma ^{-1}}\nu \sigma _{X}^{-1}(u))=$

\noindent $l\sigma _{Y}\nu \sigma _{X}^{-1}(u)$,\textit{\ i.e.}, $\sigma
_{Y}\nu \,\sigma _{X}^{-1}(lu)=l\sigma _{Y}\nu \sigma _{X}^{-1}(u)$.
Therefore, $\sigma _{Y}\nu \,\sigma _{X}^{-1}$ is a $U(L)$-module
homomorphism.\textit{\ \ \ \ \ \ }$_{\square }\medskip $

Let $\hat{\sigma}$ denote an important type of semi-inner automorphism of
the category $_{U(L)}\mathcal{M}^{0}$ of the following form:\smallskip

1) $\hat{\sigma}(F_{X})=F_{X},\;F_{X}\in |_{U(L)}\mathcal{M}^{0}|$%
,\thinspace \textit{\ i.e.}, $\hat{\sigma}$ does not change objects of the
category $_{U(L)}\mathcal{M}^{0}$;

\medskip 2) $\hat{\sigma}(\nu )=\sigma _{Y}\nu \,\sigma
_{X}^{-1}:F_{X}\longrightarrow F_{Y}$ \ for all $\nu :F_{X}\longrightarrow
F_{X}$, and $F_{X},\,F_{Y}\in $ $|_{U(L)}\mathcal{M}^{0}|$.\smallskip

In this connection, it is appropriate to mention the following.\medskip

\noindent \textbf{Remark 4.14.} Let $K$ be a commutative ring and $\sigma $
a ring automorphism of $U(L)$ such that $\sigma (K)\neq K$. Then $\hat{\sigma%
}$ is an automorphism of the category $_{U(L)}\mathcal{M}^{0}$ but it is not
a semi-inner one. Since $U(L)$ is an IBN-ring, by Theorem 3.9 $\hat{\sigma}$
is a semi-inner automorphism of the category $_{RU(L)}\mathcal{M}^{0}$ of
finitely generated free modules over the ring $U(L)$. We may say that such
an automorphism $\hat{\sigma}$ of the category $_{U(L)}\mathcal{M}^{0}$ is 
\textit{almost semi-inner}. And, therefore, all automorphisms of the
category $_{U(L)}\mathcal{M}^{0}$ are always semi-inner or almost
semi-inner. If $\sigma (K)=K$ for every $\sigma \in Aut\,U(L)$, then, in the
same way as in the proof of Theorem 4.7, one can show that all automorphisms
of $_{U(L)}\mathcal{M}^{0}$ are semi-inner.\medskip

The following result --- an analog of Theorem 3.15 and Corollary 3.16 in a
Lie algebra setting --- relates the group $Out$ $(_{U(L)}\mathcal{M}^{0})=$ $%
Aut$ $(_{U(L)}\mathcal{M}^{0})/Int$ $(_{U(L)}\mathcal{M}^{0})$ of outer
automorphisms of the category $_{U(L)}\mathcal{M}^{0}$ with the group $%
Aut\,U(L)$ of ring automorphisms of $U(L)$. \medskip

\noindent \textbf{Theorem 4.15.} \textit{Let} $L$ \textit{be a Lie algebra
over an integral domain }$K$\textit{.} \textit{Then, } $Out$ $(_{U(L)}%
\mathcal{M}^{0})=$ $Aut$ $(_{U(L)}\mathcal{M}^{0})/Int$ $(_{U(L)}\mathcal{M}%
^{0})\cong Aut\,U(L)$.\textit{\medskip }

\noindent \textbf{Proof.} By Proposition 4.2 $_{U(L)}\mathcal{M}^{0}$ is an
IBN-variety, therefore, by Proposition 2.9 and 2.7, $Aut\;(_{U(L)}\mathcal{M}%
^{0})=StAut\;(_{U(L)}\mathcal{M}^{0})Int\;(_{U(L)}\mathcal{M}^{0})$. Thus,
we have $Aut\;(_{U(L)}\mathcal{M}^{0})/Int\;(_{U(L)}\mathcal{M}^{0})\cong $

\noindent $StAut(_{U(L)}\mathcal{M}^{0})Int(_{U(L)}\mathcal{M}^{0})/Int$ $%
(_{U(L)}\mathcal{M}^{0})\cong $

\noindent $StAut$\textit{\ }$(_{U(L)}\mathcal{M}^{0})/(StAut\mathit{\ }%
(_{U(L)}\mathcal{M}^{0})\cap Int$ $(_{U(L)}\mathcal{M}^{0}))$.

We define a group homomorphism $\pi :StAut\;(_{U(L)}\mathcal{M}%
^{0})\rightarrow Aut\,U(L)$ as $\pi (\varphi )=\sigma $, where $\sigma $ is
the automorphism of $U(L)$ corresponding to $\varphi \in StAut\;(_{U(L)}%
\mathcal{M}^{0})$. Since $\varphi (\hat{\sigma})=\sigma $ for all $\sigma
\in Aut\,U(L)$, the homomorphism $\pi $ is surjective. So, we need only to
show that $Ker\,\pi =StAut\mathit{\ }(_{U(L)L}\mathcal{M}^{0})\cap Int$ $%
(_{U(L)}\mathcal{M}^{0})$.

Let $\pi (\varphi )=\sigma $, $\varphi \in StAut\mathit{\ }(_{U(L)}\mathcal{M%
}^{0})$. We will show that $\hat{\sigma}$ and $\varphi $ act in the same way
on the semigroup $End\,(U(L)x_{1})$.

Consider the endomorphism $\sigma _{x_{1}}\nu _{l}\sigma _{x_{1}}^{-1}$,
where, just for simplicity, we write $x_{1}$ for the singleton $\{x_{1}\}$,
and $l\in U(L)$. We have $\sigma _{x_{1}}\nu _{l}\sigma
_{x_{1}}^{-1}(x_{1})=\sigma _{x_{1}}\nu _{l}(x_{1})=\sigma _{x_{1}}(lx_{1})$

\noindent $=l^{\sigma }x_{1}=\nu _{{l}^{\sigma }}(x_{1}).$ By definition, $%
\varphi (\nu _{l})=\nu _{{l}^{\sigma }}$. Thus, $\hat{\sigma}=\varphi $ on $%
End\,(U(L)x_{1})$, and $\varphi _{1}=\hat{\sigma}^{-1}\varphi $ is an
automorphism of $_{U(L)}\mathcal{M}^{0}$ acting identically on the semigroup 
$End\,(U(L)x_{1})$.

Recall that the variety $_{U(L)}\mathcal{M}$ is generated by the cyclic
module $U(L)x_{1}$. Consider a homomorphism $\nu _{1}:U(L)x_{1}\rightarrow
U(L)x_{1}$ such that $\nu _{1}(x_{1})=x_{1}$. It is clear that $\varphi
_{1}(\nu _{1})=\nu _{1}$. Hence, all conditions of Theorem 5.3 (see
Appendix) hold, and therefore the automorphism $\varphi _{1}=\hat{\sigma}%
^{-1}\varphi $ is an inner automorphism of the category $_{U(L)}\mathcal{M}%
^{0}$. Thus, $\varphi =\hat{\sigma}\varphi _{1}$ is semi-inner.

Now, if $\sigma =1$ then $\varphi =\varphi _{1}$ is an inner automorphism of 
$_{U(L)}\mathcal{M}^{0}$, \textit{i.e.}, $\varphi \in StAut\mathit{\ }%
(_{U(L)}\mathcal{M}^{0})\cap Int$ $(_{U(L)}\mathcal{M}^{0})$.

Conversely, let $\varphi \in StAut\mathit{\ }(_{U(L)}\mathcal{M}^{0})\cap
Int $ $(_{U(L)}\mathcal{M}^{0})$. Consider

\noindent $\varphi _{U(L)x_{1}}$, the restriction of the morphism $\varphi $
to the cyclic module $U(L)x_{1}$. If $\nu :U(L)x_{1}\rightarrow U(L)x_{1}$,
then $\varphi _{U(L)x_{1}}(\nu )=s_{U(L)x_{1}}^{-1}\nu s_{U(L)x_{1}}$, where 
$s_{U(L)x_{1}}$ is an automorphism of $U(L)x_{1}$. Since, by Lemma 4.11,

\noindent\ $Aut\;(U(L)x_{1})=Aut\,(Kx_{1})$, we have $\varphi
_{U(L)x_{1}}(\nu )=\nu $, \textit{i.e.} $\sigma $ corresponding to $\varphi
_{U(L)x_{1}}$ is equal to $1$. This ends the proof of the theorem. \textit{\
\ \ \ \ \ }$_{\square }\medskip $

Because of the abovementioned isomorphism between the categories $_{L}%
\mathcal{M}$ and $_{U(L)}\mathcal{M}$, it is clear that $Aut(_{L}\mathcal{M}%
^{0})=Aut(_{U(L)}\mathcal{M}^{0})$ and

\noindent $Int(_{L}\mathcal{M}^{0})=Int(_{U(L)}\mathcal{M}^{0})$. From this
observation and Theorem 4.15, we deduce the following fact.\medskip\ 

\noindent \textbf{Corollary 4.16. }\textit{Let} $L$ \textit{be a Lie algebra
over an integral domain }$K$\textit{.} \textit{Then, } $Out$ $(_{L}\mathcal{M%
}^{0})=$ $Aut$ $(_{L}\mathcal{M}^{0})/Int$ $(_{L}\mathcal{M}^{0})\cong
Aut\,U(L)$.\textit{\ \ \ \ \ \ }$_{\square }\medskip $

Since for a Lie algebra $L$ over an integral domain $K$ every automorphism $%
\varphi \in Aut(_{L}\mathcal{M}^{0})$ is semi-inner, from Definition 4.6 it
follows that for every automorphism $\varphi \in Aut(_{L}\mathcal{M}^{0})$
there exists a ring automorphism $\tilde{\varphi}:U(L)\rightarrow U(L)$
corresponding to $\varphi $ (see Definition 4.6). Then, denoting by $G$ and $%
Aut_{s}\,L$ a subgroup of $Aut(_{L}\mathcal{M}^{0})$, consisting of all $%
\varphi \in Aut(_{L}\mathcal{M}^{0})$ with $\tilde{\varphi}(L)=L$, and the
group $Aut_{s}\,L$ of all semi-automorphisms of $L$, respectively, we also
have the following result.\medskip

\noindent \textbf{Proposition 4.17.} \textit{For any Lie algebra} $L$ 
\textit{over an integral domain }$K$\textit{,}

\noindent $G/Int(_{L}\mathcal{M}^{0})\cong Aut_{s}\,L$. \medskip

\noindent \textbf{Proof.} By Theorem 4.15, $G/Int(_{U(L)}\mathcal{M}%
^{0})\cong Aut_{1}\,U(L)$, where

\noindent $Aut_{1}\,U(L)$ consists of all ring automorphisms $\phi $ of $%
U(L) $ for which $\phi (L)=L$. We prove that $Aut_{1}\,U(L)=Aut_{s}\,L$. It
is sufficient to prove that every semi-automorphism $\theta =(\chi ,\omega
):L\rightarrow L$, where $\chi \in Aut\,K$ and $\omega :L\rightarrow L$ is a
ring automorphism of $L$, can be lifted up to a ring automorphism $\theta
^{\prime }:U(L)\rightarrow U(L)$. We can define a semi-automorphism $\bar{%
\chi}:L\rightarrow L$ of the algebra $L$: $\bar{\chi}\sum_{i}\alpha
_{i}e_{i}=\sum_{i}\alpha _{i}^{\chi }e_{i}$ for any $\alpha _{i}\in
K,\;e_{i}\in E$. The mapping $\chi _{1}=\bar{\chi}^{-1}\theta $ is an
automorphism of the Lie algebra $L$. It can be lifted up to automorphism $%
\theta _{1}^{\prime }$ of the algebra $U(L)$ \cite[Corollary 5.2.]%
{jacobson:jac}. Therefore, $\theta ^{\prime }=\bar{\chi}\theta _{1}^{\prime
} $ is an automorphism of the ring $U(L)$ such that a restriction $\theta
^{\prime }$ on $L$ coincides with the semi-automorphism $\theta $ of $L$. 
\textit{\ \ \ \ \ \ }$_{\square }\medskip $\medskip

\noindent \textbf{Remark 4.18.} For a $p$-Lie algebra $G$ over an integral
domain $\Phi $, it can be shown in\ the same fashion that analogous results
are valid for the categories $_{G}\mathcal{M}_{p}$ and $_{U_{p}(G)}\mathcal{M%
}_{p}$\ of free $p$-Lie modules over a $p$-Lie algebra $G$ and free modules
over $U_{p}(G)$, respectively. \textit{\ \ \ \ \ \ }$_{\square }\medskip $

In all our considerations in this section, Theorems 4.1, 4.9 and the PBWT
played quite essential role. In light of this and Remark 4.14, we conclude
this section by posting the following problem.\medskip\ 

\noindent \textbf{Problem 5.} For any Lie algebra $L$ which is not a free $K$%
-module, does there exist a non-almost-semi-inner automorphism $\varphi \in
Aut(_{L}\mathcal{M}^{0})$?\textit{\ }(Our conjecture is that the answer is
``Yes.'')\textit{\ \ \ }

\section{Appendix}

In this section, we present new and shorter proofs of several important
results from \cite{mashplts:aocfa} and \cite{masplotbe:acfla}, which are
actively used in the field and which have also been used in the present
paper. In particular, we give a categorical-algebraic proof of a
generalization of the reduction theorem (see also \cite{zhit:aganport}).

\medskip Recall \cite[Definition 10.5.1]{schu:cat} that an object $A_{0}$ ($%
A^{0}$) of a category $\mathcal{C}$ is called a \textit{separator} (\textit{%
coseparator}) if for any two different morphisms $f,g\in $ $Mor_{\mathcal{C}%
} $\textit{\ }$(A,B)$ there is a morphism $h:A_{0}\longrightarrow $ $A\ $($%
h:B\longrightarrow A^{0}$) with $fh\neq gh$ ($hf\neq hg$).\medskip

\noindent \textbf{Proposition 5.1. }(\textit{cf.} \cite[Theorem 2]%
{zhit:aganport}) \textit{Let } $A_{0}$, $A^{0}$ $\in |\mathcal{C}|$, $A_{0}$ 
\textit{be a separator, and} $\varphi \in StAut$ $(\mathcal{C})$ \textit{a
stable automorphism} \textit{such that} $\varphi (f)=f$ \ \textit{for any} $%
f\in $ $Mor_{\mathcal{C}}$\textit{\ }$(A_{0},A^{0})$\textit{. Also, suppose
that if for a bijection} $s:$ $Mor_{\mathcal{C}}$\textit{\ }$%
(A_{0},A)\longrightarrow Mor_{\mathcal{C}}$\textit{\ }$(A_{0},A)$, $A\in |%
\mathcal{C}|$ \textit{and any} $g:A\longrightarrow A^{0}$\textit{\ there
exists }$h_{g}:A\longrightarrow A^{0}$ \textit{such that} $Mor_{\mathcal{C}}$%
\textit{\ }$(A_{0},g)\circ s=Mor_{\mathcal{C}}$\textit{\ }$(A_{0},h_{g})$, 
\textit{\ then} $s=Mor_{\mathcal{C}}$\textit{\ }$(A_{0},t):Mor_{\mathcal{C}}$%
\textit{\ }$(A_{0},A)\longrightarrow Mor_{\mathcal{C}}$\textit{\ }$(A_{0},A)$
\textit{for some isomorphism} $t_{A}:A\longrightarrow A$\textit{. Then the
automorphism} $\varphi $ \textit{is inner, i.e.,} $\varphi \in Int$ $(%
\mathcal{C})$\textit{.\medskip }

\noindent \textbf{Proof}. For any $A\in |\mathcal{C}|$, the automorphism $%
\varphi $ defines the corresponding bijections $\varphi _{A}:Mor_{\mathcal{C}%
}$\textit{\ }$(A_{0},A)\longrightarrow Mor_{\mathcal{C}}$\textit{\ }$%
(A_{0},A)$ with $\varphi _{A^{0}}=Mor_{\mathcal{C}}$\textit{\ }$%
(A_{0},1_{A^{0}})$. Therefore, for any $u\in $ $Mor_{\mathcal{C}}$\textit{\ }%
$(A_{0},A)$ and $g:A\longrightarrow B$ we have $Mor_{\mathcal{C}%
}(A_{0},\varphi (g))(u)=\varphi (g)u=\varphi (g\varphi _{A}^{-1}(u))=$

\noindent $\varphi (Mor_{\mathcal{C}}\mathit{\ }(A_{0},g)(\varphi
_{A}^{-1}(u)))=\varphi _{B}(Mor_{\mathcal{C}}\mathit{\ }(A_{0},g)(\varphi
_{A}^{-1}(u))):$

\noindent $Mor_{\mathcal{C}}$\textit{\ }$(A_{0},A)\longrightarrow Mor_{%
\mathcal{C}}$\textit{\ }$(A_{0},B)$.

In particular, when $B=A^{0}$ we have $Mor_{\mathcal{C}}(A_{0},\varphi
(g))(u)=$

\noindent $Mor_{\mathcal{C}}\mathit{\ }(A_{0},g)(\varphi _{A}^{-1}(u))$.
Hence, for any $A\in |\mathcal{C}|$, there exists an isomorphism $%
t_{A}^{-1}:A\longrightarrow A$ such that $\varphi _{A}^{-1}=Mor_{\mathcal{C}%
} $\textit{\ }$(A_{0},t_{A}^{-1}):Mor_{\mathcal{C}}$\textit{\ }$%
(A_{0},A)\longrightarrow Mor_{\mathcal{C}}$\textit{\ }$(A_{0},A)$ and,
hence, $\varphi _{A}=Mor_{\mathcal{C}}$\textit{\ }$(A_{0},t_{A})$.

Thus, for any $u\in $ $Mor_{\mathcal{C}}$\textit{\ }$(A_{0},A)$ and $%
g:A\longrightarrow B$ we have $Mor_{\mathcal{C}}(A_{0},$

\noindent $\varphi (g))(u)=\varphi (g)u=\varphi (g\varphi
_{A}^{-1}(u))=\varphi (gt_{A}^{-1}u)$ $=\varphi
_{B}(gt_{A}^{-1}u)=t_{B}gt_{A}^{-1}u$; and since\ $A_{0}$ is a separator, $%
\varphi (g)=t_{B}gt_{A}^{-1}$.\textit{\ \ \ \ \ \ }$_{\square }\medskip $

Now let $\Theta $ be a variety of universal algebras and $F_{n}\in |\Theta
^{0}|$ a free algebra for some $n\in $ $\mathbf{N}$. Without loss of
generality, we may consider $F_{n}$ to be a coproduct of the $n$-copies of $%
F_{1}$, \textit{i.e.}, $F_{n}=$ $\amalg _{i=1}^{n}F_{i}\overset{\mu _{i}}{%
\leftarrow }F_{i}$, where $F_{i}=$ $F_{1}$ with the canonical injections\ $%
\mu _{i}$. Thus, there exists the codiagonal morphism $\nu _{0}\overset{def}{%
=}\amalg _{i=1}^{n}1_{F_{i}}:\amalg _{i=1}^{n}F_{i}\longrightarrow F_{1}$
such that $\nu _{0}\mu _{i}=1_{F_{i}}$ for any $i=1,\ldots ,n$. The
following observation is obvious and useful.\medskip

\noindent \textbf{Lemma 5.2.} \textit{Let} $\varphi \in StAut$\textit{\ }$%
(\Theta ^{0})$\textit{. Then,} $\varphi $ \textit{acts trivially on the
monoid}\noindent\ $Mor_{\Theta ^{0}}(F_{n},F_{n})$ \textit{and} $\varphi
(\nu _{0})=\nu _{0}$ \textit{iff} $\varphi $ \textit{acts trivially on} 

\noindent $Mor_{\Theta ^{0}}(F_{1},F_{n})$\textit{.} \medskip 

\noindent \textbf{Proof}. The result follows immediately as $\varphi $
preserves coproducts, $\nu _{0}\mu _{i}=1_{F_{1}}$, and $\nu _{0}$ is epi.%
\textit{\ \ \ \ \ \ }$_{\square }\medskip $

\noindent \textbf{Theorem 5.3. }(Reduction Theorem \cite[Theorem 3.11]%
{mashplts:aocfa} and/or \cite[Theorem 3]{masplotbe:acfla}, also cf. %
\cite[Theorem 3]{zhit:aganport}) \textit{Let a free algebra }$F_{n}$ \textit{%
generate a variety} $\Theta $\textit{, and} $\varphi \in StAut$\textit{\ }$%
(\Theta ^{0})$\textit{. If }$\varphi $ \textit{acts trivially on the monoid} 
$Mor_{\Theta ^{0}}(F_{n},F_{n})$ \textit{and} $\varphi (\nu _{0})=\nu _{0}$%
\textit{, then} $\varphi $ \textit{is inner, i.e.,} $\varphi \in Int$ $%
(\Theta ^{0})$\textit{.\medskip }

\noindent \textbf{Proof}. First, it is clear that $F_{1}$ is a separator in $%
\Theta ^{0}$. Then, suppose that $s:$ $Mor_{_{\Theta
^{0}}}(F_{1},F_{X})\longrightarrow Mor_{_{\Theta ^{0}}}(F_{1},F_{X})$, $%
F_{X}\in |\Theta ^{0}|$, is a bijection such that for any homomorphism $%
g:F_{X}\longrightarrow F_{n}$ there exists a homomorphism\textit{\ }$%
h_{g}:F_{X}\longrightarrow F_{n}$ and $Mor_{_{\Theta ^{0}}}(F_{1},g)\circ
s=Mor_{_{\Theta ^{0}}}(F_{1},h_{g})$. The bijection $s$ induces the obvious
bijection $\overline{s}:F_{X}\longrightarrow F_{X}$ on the universe $F_{X}$.
We shall show that $\overline{s}$ is an automorphism of $F_{X}$.

Indeed, let $\omega $ be a $m$-ary fundamental operation, and $a_{1},\ldots
,a_{m}\in $ $F_{X}$. Then we have two elements $\overline{s}(\omega
(a_{1},\ldots ,a_{m}))$ and $\omega (\overline{s}(a_{1}),\ldots ,$

\noindent $\overline{s}(a_{m}))$ in $F_{X}$, and, hence, for any $%
g:F_{X}\longrightarrow F_{n}$ and a corresponding $h_{g}:F_{X}%
\longrightarrow F_{n}$ we have $g(\overline{s}(\omega (a_{1},\ldots
,a_{m})))=h_{g}(\omega (a_{1},\ldots ,a_{m}))=$

\noindent $\omega (h_{g}(a_{1}),\ldots ,h_{g}(a_{m}))=\omega (g(\overline{s}%
(a_{1})),\ldots ,g(\overline{s}(a_{m})))=$

\noindent $g(\omega (\overline{s}(a_{1}),\ldots ,\overline{s}(a_{m})))$.
Thus, for any $g:F_{X}\longrightarrow F_{n}$ we always have $g(\overline{s}%
(\omega (a_{1},\ldots ,a_{m})))=g(\omega (\overline{s}(a_{1}),\ldots ,%
\overline{s}(a_{m})))$. As $F_{n}$ generates $\Theta $, by, for example, %
\cite[Theorem 13.2]{mal:as}, $F_{X}$ is a homomorphic image of a subalgebra
of a product of copies of $F_{n}$. From this and using that a free algebra\ $%
F_{X}$ is a projective object in $\Theta $, one may easily see that $F_{X}$
is a subalgebra of a product of copies of $F_{n}$ and, hence, $F_{n}$ is a
coseparator in $\Theta ^{0}$. From the latter, one immediately obtains that $%
\overline{s}(\omega (a_{1},\ldots ,a_{m}))=\omega (\overline{s}%
(a_{1}),\ldots ,\overline{s}(a_{m}))$, $\overline{s}:F_{X}\longrightarrow
F_{X}$ is an automorphism, and $s=$ $Mor_{_{\Theta ^{0}}}(F_{1},\overline{s}%
) $. Then, using Proposition 5.1 and Lemma 5.2, we end the proof.\textit{\ \
\ \ \ \ }$_{\square }\medskip $

A particular case of Theorem 5.3 is \medskip

\noindent \textbf{Theorem 5.4. }(\cite[Theorem 3.10]{mashplts:aocfa}) 
\textit{If a variety} $\Theta $ \textit{is generated by} $F_{1}$, $\varphi
\in StAut$\textit{\ }$(\Theta ^{0})$\textit{, and }$\varphi $ \textit{acts
trivially on the monoid} $Mor_{\Theta ^{0}}(F_{1},F_{1})$\textit{, then} $%
\varphi \in Int$ $(\Theta ^{0})$\textit{.\ \ \ \ \ \ }$_{\square }\medskip $

\end{document}